\documentclass[MSNbibl,number,citesort,dvips]{arxbj}

\aid{0}
\volume{19}
\issue{5A}
\pubyear{2013}
\firstpage{1880}
\lastpage{1919}
\doi{10.3150/12-BEJ434} %kopijuoti is PTS

\makeatletter
\newcommand{\reel}{{\mathbb R}}
\newcommand{\entier}{{\mathbb Z}}
\newcommand{\eps}{\varepsilon}
\newcommand{\tribuf}{{\mathcal F}}
\newcommand{\proba}{{\mathbb P}}
\newcommand{\esp}{{\mathbb E}}
\newcommand{\abar}{{\overline a}}
\newcommand{\Ntilde}{{\widetilde N}}
\newcommand{\lieal}{{\mathfrak g}}
\newcommand{\lieh}{{\mathfrak h}}
\newcommand{\liep}{{\mathfrak p}}
\newcommand{\liem}{{\mathfrak m}}
\newcommand{\atilde}{{\widetilde a}}
\newcommand{\btilde}{{\widetilde b}}
\newcommand{\Ad}{{\mathrm{Ad}}}
\newcommand{\hyperb}{{\mathbb H}}
\newcommand{\measure}{{\mathcal M}}
\newcommand{\haar}{{\mathcal H}}
\newcommand{\gener}{{\mathcal L}}
\newcommand{\generbar}{{\overline\gener}}
\newcommand{\genertilde}{{\widetilde\gener}}
\newcommand{\Xtilde}{{\widetilde X}}
\newcommand{\am}{{\mathcal A}}
\newcommand{\saut}{{\mathcal T}}
\newcommand{\pbar}{{\overline p}}
\newcommand{\ptilde}{{\widetilde p}}
\newcommand{\esptilde}{{\widetilde\esp}}
\newcommand{\Lambdatilde}{{\widetilde\Lambda}}

\newproclaim{assumption}{Assumption}
\newtheorem{theorem}{Theorem}
\newtheorem{lemma}{Lemma}

\newcommand{\SO}{\operatorname{SO}}
\newcommand{\ess}{\operatorname{ess}}
\newcommand{\eqref}[1]{(\ref{#1})}
\makeatother

\begin{document}
\begin{frontmatter}

\title{Smoothness of the law of manifold-valued Markov processes with jumps}
\runtitle{Smoothness of the law of manifold-valued processes with jumps}

\begin{aug}
%%%% inicialai - be tarpu
\author[1,2]{\fnms{Jean} \snm{Picard}\corref{}\thanksref{1,2,e1}\ead[label=e1,mark]{Jean.Picard@math.univ-bpclermont.fr}\ead[label=u1,url]{http://math.univ-bpclermont.fr/\textasciitilde picard/}} \and
\author[1,2]{\fnms{Catherine} \snm{Savona}\thanksref{1,2,e2}\ead[label=e2,mark]{Catherine.Savona@math.univ-bpclermont.fr}}
\runauthor{J. Picard and C. Savona} %% auto
\address[1]{Clermont Universit\'{e}, Universit\'{e} Blaise Pascal,
Laboratoire de Math\'{e}matiques, BP 10448,
\mbox{F-63000} Clermont-Ferrand, France. \printead{e1};\\ \printead*{e2}; \printead{u1}}
\address[2]{CNRS, UMR 6620, LM, BP 80026, F-63171 Aubi\`{e}re, France}
\end{aug}

% HISTORY:
\received{\smonth{7} \syear{2011}}

% ABSTRACT
%
\begin{abstract}
Consider on a manifold the solution $X$ of a stochastic differential
equation driven by a L\'{e}vy process without Brownian
part. Sufficient conditions for the smoothness of the law of $X_t$ are
given, with particular emphasis on noncompact
manifolds. The result is deduced from the case of affine spaces by
means of a localisation technique. The particular
cases of Lie groups and homogeneous spaces are discussed.
\end{abstract}


\end{frontmatter}

%s1 #&#
\section{Introduction}

Consider a $\reel^m$-valued L\'{e}vy process $\Lambda_t$ without
Brownian part, a $d$-dimensional manifold $M$, and the
$M$-valued solution $X_t$ of an equation
%
%e1 #&#
\begin{equation}\label{eqx}
 X_{t+\mathrm{d}t}=a(X_t,\mathrm{d}\Lambda_t)+b(X_t)\,\mathrm{d}t,\qquad
X_0=x_0
\end{equation}
for coefficients $a$ and $b$ such that $a(x,0)=x$. The precise meaning
of this equation will be given later, as well as
conditions implying the existence and uniqueness of a solution. The aim
of this article is to give sufficient
conditions ensuring the smoothness of the law of $X_t$ at any time
$t>0$. We also study more precisely the case of L\'{e}vy
processes on Lie groups, and of some classes of processes on
homogeneous spaces.

Proving the smoothness of the law of a random variable has motivated,
in the case of continuous diffusions, the
introduction of Malliavin's calculus. When $M=\reel^d$, Bismut's
approach to this calculus has proved to be useful for
processes with jumps which are solutions of equations of type \eqref{eqx}, and this topic has been intensively studied
since \cite{bismut83} and \cite{bichgravjac}. Different techniques,
each of them having its own domain of
applicability, have been introduced afterwards. These techniques can be
roughly divided into two classes.

The first class relies, as in \cite{bismut83}, on some infinitesimal
perturbations (in space or in time) on the jumps
of the L\'{e}vy process $\Lambda$. A differential calculus can be
based on these perturbations, and the associated
integration by parts formula enables to study the smoothness of the law
of $X_t$.

The second class of techniques (also when $M=\reel^d$) has been worked
out in \cite{picard96}; rather than a
differential calculus, one uses a finite difference calculus consisting
in appending and removing jumps. This is not a
differential calculus so there is no integration by parts formula in
the usual sense, but there is still a duality
formula which can be written on the Poisson space of jumps, and this
formula can be interpreted as a duality between
appending and removing jumps, see \cite{nualartvives95,picard96a};
this calculus has been applied to the smoothness of
the law of $X_t$ by \cite{picard96}. Its advantage is that no
smoothness is required for the L\'{e}vy measure of $\Lambda$;
in particular this measure may have a countable support. The proof of
the smoothness of the law is based on an
estimation of the characteristic function (Fourier transform) of $X_t$;
if this function is proved to decrease rapidly
at infinity, then $X_t$ has a $C^\infty$ density. In \cite{picard96}
(as it will be the case in this article), it is
assumed that $\Lambda$ has no Brownian part and that the coefficient
$a$ satisfies a nondegeneracy condition similar
to the ellipticity condition for continuous diffusions; notice that the
class of so-called canonical equations with a
Brownian part was studied with this finite difference calculus by \cite
{ishikkunita06}, and H\"{o}rmander type conditions
were considered by \cite{kunita09}, also for canonical equations.

The extension of \cite{picard96} to the case where $M$ is a manifold
is not immediate. The basic tool of the case
$M=\reel^d$ was the Fourier transform. When $M$ is a symmetric space,
the Fourier transform can be replaced by the
so-called spherical transform, and this has been used by \cite
{liaowang07} in order to prove the smoothness of the law
for a class of processes (but the L\'{e}vy measure was not allowed to
be purely atomic in the case without Brownian part).
However, adapting the Fourier transform technique of $\reel^d$ seems
difficult for general manifolds. When $M$ is
diffeomorphic to $\reel^d$, we can of course apply the result for
$\reel^d$, but the assumptions needed for this result
are generally not translated into canonical assumptions on $M$: they
depend on the choice of the diffeomorphism, and
they may be hard to verify, see in Section \ref{isotrop} the example
of the hyperbolic space. Actually, even the case
where $M$ is an open subset of $\reel^d$ is not trivial.

In order to get around these difficulties, our aim here is to apply
some localisation techniques; we consider an atlas
of $M$ and use the results of the affine case on each local chart.
However, it is known that localisation is made
difficult by the presence of jumps. In \cite{picsav00}, we have
applied such a technique in order to study the
smoothness of harmonic functions on some domain $D$ of $\reel^d$, and
it appears that these functions are not always
$C^\infty$, even for processes with a $C^\infty$ density; their order
of regularity depends on the number of jumps
needed to exit $D$, so that the smaller the jumps are, the smoother the
function is. Our plan is therefore, first, to
apply the localisation technique and obtain the $C^\ell$ regularity of
the density when big jumps are removed, and then
check that adding big jumps does not destroy this smoothness. In order
to conclude, when the manifold is not compact,
we have to make an assumption on the size of these big jumps, namely
that for any relatively compact subset $U$ of $M$,
the set of points from which the process can jump into $U$ is
relatively compact (one cannot jump from a very distant
location). This condition can be viewed as dual to the condition of
\cite{picsav00} for the smoothness of harmonic
functions. The result (Theorem \ref{mainth}) is given in Section \ref
{main}, after the particular case of canonical
equations (Theorem \ref{canth}); it is proved in Section \ref{proofs}.

Then we relax this ``big jumps'' condition in three cases:
\begin{itemize}
\item when the L\'{e}vy kernel for these big jumps is smooth (Theorem
\ref{smoothj}),
\item when $M$ is a Lie group $G$ and $X_t$ is a L\'{e}vy process on it,
\item when $M$ is an homogeneous space $G/H$ and $X_t$ is obtained by
projecting on $M$ a L\'{e}vy process on $G$.
\end{itemize}
Lie groups and homogeneous spaces are the purpose of Section \ref{liehom}.

In Section \ref{exam}, we give some examples, and also some
counterexamples where the ``big jumps'' condition is not
satisfied, and the smoothness of the law fails.

%s2 #&#
\section{The smoothness result on manifolds}\label{main}

In this section, we give a precise meaning for equation \eqref{eqx},
but before that we introduce the manifold $M$ and
the L\'{e}vy process $\Lambda_t$. Then we give the assumptions on the
equation, give the smoothness results about the law
of the solution, and make some comments on our model. Actually, before
explaining the general case, we state the main
result in the particular case of so-called canonical equations;
assumptions are indeed much easier to write in this
case. The proofs are postponed to Section \ref{proofs}.

%s2.1 #&#
\subsection{The manifold}\label{manif}

The manifold $M$ is supposed to be Hausdorff, separable, paracompact,
$C^\infty$ and of dimension $d$; we do not
suppose that it is connected; it may have an at most countable number
of connected components, and our processes will
be allowed to jump from a component to another. The tangent bundle is
denoted by $TM=\bigcup_xT_xM$. If $M$ is not
compact, we consider its one-point compactification $M\cup\{\infty\}
$; if $M$ is compact, $\infty$ is a point which is
disconnected from $M$. This additional point will be viewed as a
cemetery point; this means that real functions $f$ on
$M$ are extended by putting $f(\infty)=0$, and that a process on
$M\cup\{\infty\}$ which hits $\infty$ or tends to
$\infty$ stays at that point forever.

Under these assumptions, one can embed $M$ into an affine space $\reel^N$ by means of Whitney's theorem, and one can
consider Riemannian metrics on $M$ by using a partition of unity;
Riemannian distances are of course only defined
between points in the same connected component, and they otherwise take
the value $+\infty$. We would like to have
intrinsic assumptions, which do not depend on a particular embedding or
a particular metric. This will be possible
because Riemannian distances are equivalent on compact subsets, so a
metric condition which is supposed to hold
uniformly on compact subsets will actually not depend on the choice of
the metric. For instance, if $f$ is a $C^\ell$
function and $K$ is a compact subset of $M$, if we denote by $D^jf(x)$
the iterated differential of $f$ (this is a
multi-linear form on $(T_xM)^j$), we can define
\[
|f|_{K,\ell}=\sum_{j=0}^\ell
\sup_{x\in K}\bigl|D^jf(x)\bigr|
\]
if we have chosen a Riemannian metric on $M$, but changing the metric
will lead to an equivalent semi-norm; by allowing
$K$ to vary, we obtain a Fr\'{e}chet space $C^\ell(M)$. A signed
measure $\nu$ on $M$ is said to be absolutely continuous,
respectively, $C^\ell$, if its restrictions to local charts are
absolutely continuous, respectively, have $C^\ell$
densities with respect to the Lebesgue measure; it is said to have
positive density if it has (strictly) positive
density with respect to the Lebesgue measure on charts; this does not
depend on the atlas; let $\measure^{\,\ell}(M)$ be
the set of $C^\ell$ measures. If we choose a $C^\infty$ reference
measure $\mathrm{d}x$ with positive density, we can define the
family of semi-norms
\[
|\nu|_{K,\ell}=|\mathrm{d}\nu/\mathrm{d}x|_{K,\ell}.
\]
Changing the reference measure and the Riemannian metric will not
change the topology of $\measure^{\,\ell}(M)$. We also
let
\[
\measure_K^{\,\ell}=\bigl\{\nu\in\measure^{\,\ell}(M);
\nu\bigl(K^c\bigr)=0\bigr\}
\]
if $K$ is a compact subset of $M$.

The notation $U\Subset V$ for open subsets of $M$ will mean that $U$ is
relatively compact in~$V$.

%s2.2 #&#
\subsection{The L\'{e}vy process}

Let us first recall some basic facts about L\'{e}vy processes $\Lambda_t$ with values in $\reel^m$. A L\'{e}vy measure is a
measure $\mu$ on $\reel^m\setminus\{0\}$ such that $\int(|\lambda
|^2\wedge1)\mu(\mathrm{d}\lambda)<\infty$, and a L\'{e}vy process
$\Lambda_t$ without Brownian part and with L\'{e}vy measure $\mu$ is
a process which can be written by means of the
L\'{e}vy--It\^{o} representation formula
%
%e2 #&#
\begin{equation}\label{levyito}
 \Lambda_t=\kappa t+\int_0^t
\int_{\{|\lambda|\le1\}}\lambda\Ntilde (\mathrm{d}s,\mathrm{d}\lambda) +\int
_0^t\int_{\{|\lambda|>1\}}\lambda N(\mathrm{d}s,\mathrm{d}
\lambda)
\end{equation}
for some $\kappa\in\reel^m$, where the random measure $N(\mathrm{d}t,\mathrm{d}\lambda
)=\sum_t\delta_{(t,\Delta\Lambda_t)}$ is a Poisson
measure on $\reel_+\times\reel^m$ with intensity $\mathrm{d}t \mu(\mathrm{d}\lambda
)$, and
$\Ntilde(\mathrm{d}t,\mathrm{d}\lambda)=N(\mathrm{d}t,\mathrm{d}\lambda)-\mathrm{d}t \mu(\mathrm{d}\lambda)$ is the
compensated Poisson measure. When
$\int(|\lambda|\wedge1)\mu(\mathrm{d}\lambda)<\infty$, then $\Lambda_t$
has finite variation, and \eqref{levyito} can be
simplified as
\[
\Lambda_t=\kappa_0 t+\int_0^t
\int_{\reel^m}\lambda N(\mathrm{d}s,\mathrm{d}\lambda).
\]
The relation between $\kappa$ and $\kappa_0$ is easily written, and
$\Lambda_t$ is a pure jump process when
$\kappa_0=0$. We will assume an approximate self-similarity condition
and a non-degeneracy condition on the L\'{e}vy
process written as follows.

%as1 #&#
\begin{assumption}\label{aslevy}
There exist some $0<\alpha<2$ and some positive $c$ and $C$ such that
\[
c\rho^{2-\alpha}|u|^2\le\int_{\{|\lambda|\le\rho\}}\langle
\lambda,u\rangle^2\mu(\mathrm{d}\lambda)\le C\rho^{2-\alpha}|u|^2
\]
for $u\in\reel^m$, $0<\rho\le1$, and where $\langle \cdot,\cdot\rangle$ is
the Euclidean inner product. If $\alpha=1$, we
suppose moreover that
\[
 \operatorname{\lim\sup}\limits_{\eps\downarrow0}\biggl|\int_{\{\eps\le|\lambda|\le1\}
}\lambda\mu(\mathrm{d}\lambda)\biggr|<
\infty.
\]
If $\alpha<1$ (finite variation case), we suppose that $\Lambda_t$ is
pure jump ($\kappa_0=0$).
\end{assumption}

The process $\Lambda$ has finite variation if and only if $\alpha<1$.
Notice that no smoothness is assumed on $\mu$; it
can, for instance, have a countable support. The additional conditions
of the cases $\alpha=1$ and $\alpha<1$ are needed
to apply the results of \cite{picard97,picsav00}.

As an example, we can consider stable L\'{e}vy processes of index
$\alpha$, or more generally semi-stable processes (see
Chapter 3 of \cite{sato99} for definitions and properties concerning
these processes). If $\Lambda_t$ is a semi-stable
process and if $\mu$ is not supported by a strict subspace of $\reel^m$, then Assumption \ref{aslevy} is satisfied when
$1<\alpha<2$. When $0<\alpha\le1$, the additional conditions mean
that $\Lambda_t$ should be a strictly semi-stable
process. Actually, Assumption \ref{aslevy} is only concerned with
small jumps, so that we can also consider for
instance truncated semi-stable processes (where jumps greater than some
value are removed).

%s2.3 #&#
\subsection{The equation}

Let us now introduce the process $X_t$, solution of \eqref{eqx}. When
$M=\reel^d$, the meaning of this equation is
%
%e3 #&#
\begin{eqnarray}\label{eqeuclid}
X_t&=&X_u+\int
_u^t\abar(X_{s-})\,\mathrm{d}
\Lambda_s+\int_u^tb(X_s)\,\mathrm{d}s\nonumber
\\[-8pt]
\\[-8pt]
&& {}+\sum_{u<s\le t}\bigl(a(X_{s-},\Delta
\Lambda_s)-X_{s-}-\abar (X_{s-})\Delta
\Lambda_s\bigr)\nonumber
\end{eqnarray}
for $u\le t$, where $\abar(x)$ is the differential at 0 of $\lambda
\mapsto a(x,\lambda)$. Under convenient smoothness
conditions, the It\^{o} integral with respect to $\Lambda$ is well
defined, the sum converges, and the equation has a
unique solution for any initial condition $X_0=x_0$. Notice that the
jumps are given by
$X_t=a(X_{t-},\Delta\Lambda_t)$. It\^{o}'s formula enables to write
%
%e4 #&#
\begin{eqnarray}\label{fxt}
f(X_t)&=&f(x_0)+\int
_0^t(Df \abar) (X_{s-})\,\mathrm{d}
\Lambda_s+\int_0^t(Df b)
(X_s)\,\mathrm{d}s\nonumber
\\[-8pt]
\\[-8pt]
&&{}+\sum_{0<s\le t}\bigl((f\circ a) (X_{s-},
\Delta\Lambda_s)-f(X_{s-})-(Df\abar) (X_{s-})
\Delta\Lambda_s\bigr)\nonumber
\end{eqnarray}
for smooth functions $f$, and where $Df$ denotes the differential of
$f$.\eject

This formula can be used to give a meaning to \eqref{eqx} when $M$ is
a manifold. We consider a coefficient $a$ and a
vector field $b$,
\[
a\dvtx M\times\reel^m\to M\cup\{\infty\},\qquad b\dvtx M\to TM.
\]
We again let $\abar(x)$ be the differential of $\lambda\mapsto
a(x,\lambda)$ at 0, so that $\abar(x)$ is a linear map
from $\reel^m$ into $T_xM$:
%
%e5 #&#
\begin{equation}\label{abar}
 \abar\dvtx M\to L\bigl(\reel^m,TM\bigr),\qquad x\mapsto
\partial_\lambda|_{\lambda
=0}a(x,\lambda).
\end{equation}
We say that $(X_t;t\ge0)$ is a solution of \eqref{eqx} if $X_t$ is a
c\`{a}dl\`{a}g (right continuous with left limits)
process with values in $M\cup\{\infty\}$, which is adapted with
respect to the completed filtration of $\Lambda_t$, and
such that for any smooth function $f$ and any compact subset $K$ of
$M$, equation \eqref{fxt} holds true up to the
first exit time from $K$. We say that the process dies at time $t$ if
$X_s\in M$ for $s<t$ and $X_t=\infty$. This
occurs if either $X_{t-}=\infty$ (the solution of the equation
explodes at time $t$), or if $X_{t-}\in M$ and
$X_t=a(X_{t-},\Delta\Lambda_t)=\infty$ (the process is killed at a
jump of $\Lambda$).

If we use a proper embedding ${\mathcal I}\dvtx M\to\reel^N$, we can apply
\eqref{fxt} to the components of the embedding
$\mathcal I$ and deduce the equation that should be satisfied on
${\mathcal I}(M)$ by ${\mathcal I}(X_t)$; in order to
obtain an equation on $\reel^N$, we have to extend the coefficients
out of ${\mathcal I}(M)$; it is then sufficient to
solve the equation in $\reel^N$ and prove that the solution remains in
$M$. This is this point of view which is
generally used in order to prove the existence and uniqueness of a
solution, see \cite{cohen96} in the case without
killing ($a(x,\lambda)\in M$ for any $x\in M$). We are here in a
slightly different framework (possible killing), and
shall prefer to give a proof for the existence and uniqueness by means
of local charts; this is because the proof based
on local charts is then used for studying the smoothness of the law of
$X_t$ (it is clear that the embedded process
does not have a density in $\reel^N$, so proving the existence of a
density on $M$ cannot be made through a Malliavin
calculus on $\reel^N$). The equation on a local chart can also be
deduced from \eqref{fxt}.

The infinitesimal generator of $X_t$ is
%
%e6 #&#
\begin{eqnarray}\label{infgen}
 \gener f(x)&=&Df(x) \bigl(b(x)+\abar(x)\kappa\bigr) +\int
_{\{|\lambda|\le1\}}\bigl(f\bigl(a(x,\lambda)\bigr)-f(x)-Df(x)\abar (x)
\lambda\bigr)\mu(\mathrm{d}\lambda)
\nonumber
\\[-8pt]
\\[-8pt]
&&{}+\int_{\{|\lambda|>1\}}\bigl(f\bigl(a(x,\lambda)\bigr)-f(x)\bigr)\mu(\mathrm{d}
\lambda)\nonumber
\end{eqnarray}
for $f$ bounded and $C^\infty$. In particular, if $f$ is constant on a
neighbourhood of $x$, then
\[
\gener f(x)=\int_M\bigl(f(y)-f(x)\bigr)
\mu_x(\mathrm{d}y),
\]
where $\mu_x$ is the image of the L\'{e}vy measure $\mu$ by $\lambda
\mapsto a(x,\lambda)$. This measure describes the
intensity of jumps; it is the L\'{e}vy kernel of the process $X$.

%s2.4 #&#
\subsection{Canonical equations}\label{caneq}

Up to now, we have been given an equation \eqref{eqx} on the manifold
$M$, and have explained the rigourous meaning of
this equation; in this explanation we need the function $\abar$ given
by \eqref{abar}. It is clear that many functions
$a$ are associated to the same $\abar$. However, there is a particular
class of equations, called canonical equations,
and which were introduced by \cite{marcus78} (see also \cite
{kurtzparpro95}), for which $a$ and $\abar$ are in
one-to-one correspondence.

Let us first consider a smooth field $\abar(x)\in L(\reel^m,T_xM)$,
and let
%
%e7 #&#
\begin{equation}\label{canonic}
 a(x,\lambda)=x_\lambda(1) \qquad \mbox{for } x_\lambda(t)=x+
\int_0^t\abar \bigl(x_\lambda(s)\bigr)
\lambda \, \mathrm{d}s,
\end{equation}
assuming that the solution of this ordinary differential equation does
not explode. Then it is easily seen that $a$ and
$\abar$ are related to each other by \eqref{abar}, and $x\mapsto
a(x,\lambda)$ is a diffeomorphism of $M$ onto itself
with inverse $x\mapsto a(x,-\lambda)$. Notice that $x$ and
$a(x,\lambda)$ are in the same connected component, so the
study can be reduced to connected manifolds.

For canonical equations, the assumptions needed for our main result
(Theorem \ref{mainth} below), or at least a
sufficient condition ensuring that they are satisfied, can be written
in the following simple form.

%th1 #&#
\begin{theorem}\label{canth}
Let $\Lambda$ be a L\'{e}vy process satisfying Assumption \ref
{aslevy}, and let $\abar$ and $b$ be $C^\infty$ functions on
$M$, with values, respectively, in $L(\reel^m,TM)$ and $TM$. Consider
equation \eqref{eqx} with coefficient $a$ given by
\eqref{canonic}, assuming that the ordinary differential equation
never explodes on $M$. We suppose that the jumps of
$\Lambda$ are bounded, and that the linear map $\abar(x):\reel^m\to
T_xM$ is surjective for any $x$; if $\alpha<1$ we
also suppose that $b=0$. Then \eqref{eqx} has for any initial
condition $x_0$ a unique solution $X_t$, and the law of
$X_t$ is $C^\infty$ for any $t>0$.
\end{theorem}

%s2.5 #&#
\subsection{Assumptions on the equation}

We now return to the case of a general coefficient $a$. Let us give the
assumptions on equation \eqref{eqx} which will
imply the existence and uniqueness of a solution, and the smoothness of
the law of this solution.

%as2 #&#
\begin{assumption}\label{ascoeff}
The conditions on the coefficients $a$ and $b$ are as follows.
\begin{enumerate}
\item Consider, for any $\eps>0$, the map $a_\star^\eps$ which sends
a measure $\nu$ to the measure on $M$
%
%e8 #&#
\begin{equation}\label{aeps}
 \bigl(a_\star^\eps\nu\bigr) (A)=\int\int
1_A\bigl(a(x,\lambda)\bigr)1_{\{|\lambda|>\eps\}
}\nu(\mathrm{d}x)\mu(\mathrm{d}\lambda).
\end{equation}
Let $K$ be any compact subset of $M$. If $\eps$ is small enough, then
$a_\star^\eps$ is a continuous map from
$\measure^{\,\ell}_K$ into $\measure^{\,\ell}(M)$ (see Section \ref{manif}
for the definition of these spaces).
\item Let $K$ be any compact subset of $M$. There exists $\eps>0$ such
that $x\mapsto a(x,\lambda)$ is $C^\infty$
on $K$ for $\mu$-almost any $\lambda$ such that $|\lambda|\le\eps
$.\vadjust{\goodbreak} The function $\abar$ given by \eqref{abar}
is assumed to exist and to be $C^\infty$. Letting $D$ be the
differentiation operator on $M$ with respect to
$x$, there exists some $\alpha\vee1<\gamma\le2$ such that
\[
\bigl|D^j\bigl(f\bigl(a(x,\lambda)\bigr)-f(x)-Df(x)\abar(x)\lambda\bigr)\bigr|
\le C_{f,j,K}|\lambda|^\gamma
\]
for any $C^\infty$ real function $f$, for $x\in K$, for $\mu$-almost
any $|\lambda|\le\eps$, and for any $j$.
\item The coefficient $b$ is a $C^\infty$ vector field on $M$ with
values in $TM$. In the case $\alpha<1$, we
suppose that $b=0$.
\item The linear map $\abar(x)\in L(\reel^m,T_xM)$ is surjective for
any $x$.
\end{enumerate}
\end{assumption}

Here are some comments about these four conditions.

\emph{First condition.} This condition states that jumps preserve the
smoothness of the law of the process.
Suppose that $x\mapsto a(x,\lambda)$ is a $C^\ell$ diffeomorphism of
$M$ onto itself with inverse $y\mapsto
a^{-1}(y,\lambda)$. Fix a Riemannian metric on $M$ and the associated
measure $\mathrm{d}x$; then the density $p_\star(y)$ of
$a_\star^\eps\nu$ is obtained from the density $p$ of $\nu$ by
means of the classical formula
\[
p_\star(y)=\int_{\{|\lambda|>\eps\}}\bigl|\det Da^{-1}(y,
\lambda )\bigr|p\bigl(a^{-1}(y,\lambda)\bigr)\mu(\mathrm{d}\lambda)
\]
(the determinant is computed for orthonormal bases on the tangent
spaces). It is therefore sufficient to estimate the
derivatives of $a^{-1}$; the condition is for instance satisfied for
the canonical equations of Theorem \ref{canth}.
Moreover, the process can be killed when it quits some open subset.
More precisely, consider a process on a manifold
$M_0$ associated to an equation with coefficient $a_0$, let $M$ be an
open subset of $M_0$, and kill the process when
it quits $M$. This process is obtained by considering the equation with
coefficient
%
%e9 #&#
\begin{equation}\label{killj}
 a(x,\lambda)= %
\cases{ a_0(x,\lambda) & \quad if $a_0(x,\lambda)\in M$,
\cr
\infty  &\quad otherwise. }
\end{equation}
Then $a_\star^\eps$ is continuous if the same property holds true for
$a_0$, because the map which sends a measure on
$M_0$ to its restriction to $M$ is continuous from $\measure^{\,\ell}
(M_0)$ to $\measure^{\,\ell}(M)$. We shall however notice
that Assumption \ref{asjump} below often fails for this example. Other
examples will be given in Section \ref{jcoeff}.

\emph{Second condition.} In this condition, the smoothness of
$a(x,\lambda)$ is assumed with respect to $x$ for
$\lambda$ small, but no smoothness is assumed with respect to $\lambda
$ except for $\lambda\to0$. If however $a$ is
smooth in $(x,\lambda)$, then the condition is satisfied for $\gamma
=2$. In particular, canonical equations of Section
\ref{caneq} satisfy the condition.

\emph{Third condition.} The additional assumption on $b$ in the finite
variation case $\alpha<1$ means that
the solution $X$ of our equation is a pure jump process; it is required
in order to apply the results of
\cite{picard97,picsav00}.

\emph{Fourth condition.} The surjectivity of $\abar(x)$ is a
nondegeneracy condition. This condition says
that small jumps go in all the directions, and is similar to the
ellipticity condition for continuous diffusions.

Let us now give the additional assumption concerning big jumps. It is
stating, roughly speaking, that the process
cannot come from a very distant point by jumping.

%as3 #&#
\begin{assumption}\label{asjump}
If $U$ is relatively compact, then
\[
a^{-1}(U)=\bigl\{x;\mu\bigl\{\lambda;a(x,\lambda)\in U\bigr\}>0\bigr\}
\]
is also relatively compact.
\end{assumption}

This assumption is trivially satisfied when $M$ is compact. It is also
satisfied in Theorem \ref{canth} because the
jumps of $\Lambda$ are in some ball $B$ of $\reel^m$, so
\[
a^{-1}(U)\subset a(U\times B),
\]
and the relative compactness of this set follows from the relative
compactness of $U\times B$ and the continuity of
$a$. The assumption often fails when the process is obtained by killing
as explained in \eqref{killj}; difficulties
generally arise when the original process can jump from $M_0\setminus
M$ into $M$; a~counterexample, showing that the
law of $X_t$ is not always smooth in this case, will be given in
Section \ref{killed}.

%s2.6 #&#
\subsection{The results}

The main smoothness result for the solution of \eqref{eqx} is the
following one.

%th2 #&#
\begin{theorem}\label{mainth}
Under Assumptions \ref{aslevy} and \ref{ascoeff}, the equation \eqref
{eqx} has a unique solution $X_t$ for any initial
condition $X_0=x_0$, and the law of $X_t$ is absolutely continuous for
any $t>0$. If moreover Assumption \ref{asjump}
holds true, the law of $X_t$ is $C^\infty$ for $t>0$. More precisely,
if $K$ is a compact subset of $M$, if $t_0>0$,
and if $p(t,x_0,x)$ is the density of $X_t$ with respect to some
$C^\infty$ reference measure with positive density,
then the derivatives of $p$ with respect to $x$ satisfy
%
%e10 #&#
\begin{equation}\label{djpt}
 \bigl|D^jp(t,x_0,x)\bigr|\le C_{j,K}t^{-(d+j)/\alpha}
\end{equation}
uniformly for $0<t\le t_0$, $x_0$ in $M$ and $x$ in $K$.
\end{theorem}

In the case of canonical equations, we see that Theorem \ref{canth} is
a corollary of this result. In this result, jumps
of $\Lambda$ were supposed to be bounded. When $M$ is compact but
jumps of $\Lambda$ are unbounded, we see that the
condition which may cause problems is the first condition of Assumption
\ref{ascoeff}. More precisely, if the vector
field $\abar(\cdot)\lambda$ has for some $\lambda=\lambda_1$ a stable
equilibrium $x_1$, then jumps
$\Delta\Lambda=c\lambda_1$ for large $c$ concentrate the mass near
$x_1$, so that the density of $a_\star^\eps\nu$ may
be unbounded near $x_1$.

When the big jumps condition (Assumption \ref{asjump}) is not
satisfied, the conclusion of the theorem still holds if
these big jumps are smooth; some other cases will be studied in Section
\ref{liehom}, and on the other hand,
counterexamples will be given in Section \ref{exam}.

%th3 #&#
\begin{theorem}\label{smoothj}
Under Assumptions \ref{aslevy} and \ref{ascoeff}, suppose that there
exists a decomposition $\mu=\mu^\flat+\mu^\sharp$
of the L\'{e}vy measure such that only $\mu^\flat$ satisfies
Assumption \ref{asjump}, and $\mu^\sharp$ is finite. Let
$\mu_x^\sharp$ be the image of $\mu^\sharp$ by the map $\lambda
\mapsto a(x,\lambda)$; suppose that $\mu_x^\sharp$ is
$C^\infty$, and that $x\mapsto|\mu_x^\sharp|_{K,\ell}$ is bounded
on $M$ for any compact $K$ and any $\ell$
(see the definition of this family of semi-norms in Section \ref
{manif}). Then the solution $X_t$ of \eqref{eqx} has a
$C^\infty$ law for $t>0$, satisfying \eqref{djpt}.
\end{theorem}

The kernel $\mu_x^\sharp$ is the part of the L\'{e}vy kernel $\mu_x$
coming from $\mu^\sharp$. The assumption of Theorem~\ref{smoothj} therefore requires the part of the L\'{e}vy kernel for
jumps coming from distant locations to be smooth,
whereas the L\'{e}vy kernel could be purely atomic in Theorem \ref{mainth}.

%s2.7 #&#
\subsection{About the jump coefficient}\label{jcoeff}

In this article, we study Markov processes $X_t$ with infinitesimal
generator of the form \eqref{infgen}, with
coefficients $a$ and $b$ satisfying some smoothness assumptions. This
covers a class of Markov processes, but not all
of them. In particular, in this model, the set of (times of) jumps of
$X$ is contained in the set of jumps of the
driving L\'{e}vy process $\Lambda$; the inclusion may be strict, since
one may have $a(x,\lambda)=x$ for some
$\lambda\ne0$, but introducing such a behaviour generally destroys
the smoothness of $a$. Thus, $X$ and $\Lambda$ have
more or less the same times of jumps, and the rate of jump of $X$ is
not allowed to depend on its present state $x$.
This is a drawback of this approach, as well as of other approaches
based on different versions of Malliavin's
calculus. A weak dependence of the rate of jumps with respect to $x$
can however be obtained through Girsanov
transformation, see \cite{picsav00}, and cases of more general
dependence have been studied (under frameworks which
are different from ours) by \cite{fournier08,ballyclement10,kulik11}.
We now verify that such a dependence is possible
if we drop the assumption of smoothness of $a$, and that this is
compatible with our assumptions if only finitely many
jumps are concerned by this behaviour.

Suppose that the generator of $X$ is $\gener+\gener\,'$, where $\gener
$ satisfies the assumptions of Theorem \ref{mainth}
or~\ref{smoothj}, and
\[
\gener\,'f(x)=\int_{M\cup\{\infty\}}\bigl(f(y)-f(x)\bigr)
\mu_x'(\mathrm{d}y),
\]
where $x\mapsto\mu_x'(M\cup\{\infty\})$ is finite and $C^\infty$.
We also have to assume that the kernel $\mu_x'$ is
Borel measurable, that $\nu\mapsto\nu'$ with
%
%e11 #&#
\begin{equation}\label{nuprime}
 \nu'(A)=\int\mu_x'(A)
\nu(\mathrm{d}x)
\end{equation}
is continuous from $\measure_K^{\,\ell}$ into $\measure^{\,\ell}(M)$ (as in
Assumption \ref{ascoeff}), and that jumps
corresponding to this kernel satisfy Assumption \ref{asjump} (the
extension with the assumptions of Theorem
\ref{smoothj} is also possible).

In order to prove that this case enters our framework, we apply the
fact that the measurable space $M\cup\{\infty\}$
can be viewed as a Borel subset of $\reel_+$ (it is a Lusin space, see
\cite{delmey75}); thus we can view $\mu_x'$ as
a measure on $\reel_+$, and if we define
\[
\atilde(x,u)=\inf\bigl\{v\ge0;\mu_x'\bigl([0,v]\bigr)
\ge u\bigr\}\in\reel_+\cup\{ +\infty\},
\]
the image of the Lebesgue measure on $\reel_+$ by $\atilde(x,\cdot)$ is
the measure $\mu_x'$ on $\reel_+$, plus an infinite
mass at $+\infty$; notice in particular that $\atilde(x,u)=+\infty$
if $u$ is large enough. Let us then introduce a
real symmetric L\'{e}vy process $\Lambda_t'$ with L\'{e}vy measure
$\mu'(\mathrm{d}\lambda')=\alpha|\lambda'|^{-\alpha-1}\mathrm{d}\lambda'$,
independent of $\Lambda$; then $(\Lambda,\Lambda')$ is a L\'{e}vy
process on $\reel^{m+1}$ with L\'{e}vy measure
$\mu(\mathrm{d}\lambda)\delta_0(\mathrm{d}\lambda')+\delta_0(\mathrm{d}\lambda)\mu'(\mathrm{d}\lambda')$ satisfying Assumption \ref{aslevy}. We then
consider the equation driven by $(\Lambda,\Lambda')$, with
coefficient $a(x,\lambda,\lambda')$, where $a(x,\lambda,0)$
is the coefficient associated to the part $\gener$ of the generator, and
\[
a\bigl(x,0,\lambda'\bigr)= %
\cases{
\atilde\bigl(x,{
\lambda'}^{-\alpha}\bigr)& \quad if $\lambda'>0$ and
this quantity is finite,\cr
x  &\quad otherwise.}
\]
In particular, $a(x,0,\lambda')=x$ if $\lambda'<\mu_x'(M\cup\{
\infty\})^{-1/\alpha}$. The image of $\mu'$ by $a(x,0,\cdot)$
is $\mu_x'$ plus an infinite mass at $x$, so the solution of the
equation has generator $\gener+\gener\,'$ as required.
On the other hand, if $K$ is a compact subset of $M$ and if
\[
\eps<\Bigl(\sup_{x\in K}\mu_x'\bigl(M\cup\{
\infty\}\bigr)\Bigr)^{-1/\alpha},
\]
then, with the notations \eqref{aeps} and \eqref{nuprime},
\[
\bigl(a(\cdot,0,\cdot)_\eps^\star\nu\bigr) (\mathrm{d}x)=\nu'(\mathrm{d}x)+
\bigl(\eps^{-\alpha}-\mu_x'\bigl(M\cup\{\infty\}
\bigr)\bigr)\nu(\mathrm{d}x).
\]
Thus these jumps satisfy the first part of Assumption \ref{ascoeff},
though $x\mapsto a(x,0,\lambda')$ is generally not
continuous.

We have already seen in \eqref{killj} that we can consider hard
killing of a process (we kill it when it hits an
obstacle), but this may cause difficulties with Assumption \ref
{asjump}. With the construction we have just described,
we can also consider soft killing where the process is killed at some
rate $h(x)\ge0$ depending smoothly on $x$; this
means that we add the term $\gener\,'f(x)=-h(x)f(x)$ to $\gener f(x)$,
and $\mu_x'$ is the mass $h(x)$ at $\infty$; in
this case, the measure $\nu'$ of \eqref{nuprime} is the zero measure
on $M$, so $\nu\mapsto\nu'$ is trivially
continuous.

%s3 #&#
\section{\texorpdfstring{Proof of Theorems \protect\ref{mainth} and \protect\ref{smoothj}}{Proof of Theorems 2 and 3}}\label{proofs}

The two theorems are proved in several steps.

%s3.1 #&#
\subsection{Construction of the solution}

In order to prove the existence and uniqueness of a solution of \eqref
{eqx}, we first write the equation on a local
chart $(U,\Phi)$, where $\Phi$ is a diffeomorphism from an open
subset $U$ of $M$ onto an open subset $V$ of $\reel^d$.
We can restrict ourselves to atlases such that $\Phi$ is the
restriction to $U$ of a smooth map on $M$. If $\tau$ is
the exit time of $X$ from $U$, then \eqref{fxt} applied to the
components of $\Phi$ shows that $Y_t=\Phi(X_t)$,
$t<\tau$, should be a solution of an equation \eqref{eqeuclid} on $V$
with new coefficients
%
%e12 #&#
\begin{eqnarray}\label{aphi}
a_\Phi(y,\lambda)&=&\Phi\bigl(a\bigl(
\Phi^{-1}(y),\lambda\bigr)\bigr),\qquad b_\Phi(y)=D\Phi\bigl(
\Phi^{-1}(y)\bigr)b\bigl(\Phi^{-1}(y)\bigr),\nonumber
\\[-8pt]
\\[-8pt]
\abar_\Phi(y)&=&D\Phi\bigl(\Phi^{-1}(y)\bigr)\abar\bigl(
\Phi^{-1}(y)\bigr).
\nonumber
\end{eqnarray}

More precisely, $X$ is solution of \eqref{eqx} if it is a c\`{a}dl\`
{a}g process on $M\cup\{\infty\}$ with initial condition
$X_0=x_0$, satisfying the conditions:
\begin{itemize}
\item For any local chart $(U,\Phi)$ and for any time $u$, if $\tau$
is the exit time after $u$ of $X$ from $U$,
the process $Y_t=\Phi(X_t)$ satisfies the equation \eqref{eqeuclid}
with coefficients $(a_\Phi,b_\Phi)$ on
$\{u\le t<\tau\}$.
\item The jumps of $X$ are given by $X_t=a(X_{t-},\Delta\Lambda_t)$.
\item If $X_{t-}$ or $X_t$ is at $\infty$, then $X_s=\infty$ for any
$s\ge t$ ($\infty$ is a cemetery point).
\end{itemize}
In order to solve the equation \eqref{eqx} up to the first exit time
from $U$, we shall have to extend coefficients
$(a_\Phi,b_\Phi)$ out of $V$, and solve the resulting equation \eqref
{eqeuclid} on $\reel^d$.

%le1 #&#
\begin{lemma}\label{exuniq}
The equation \eqref{eqx} has a unique solution $X_t$ for any initial
condition $x_0$.
\end{lemma}

\begin{pf}
Consider open subsets of $M$ with relatively compact inclusions
$U_1\Subset U_2\Subset U_3$. We suppose that there
exists a diffeomorphism $\Phi$ from $U_3$ onto an open subset $V_3$ of
$\reel^d$, so that $(U_3,\Phi)$ is a local
chart. We define $V_1=\Phi(U_1)$, $V_2=\Phi(U_2)$, and let $h\dvtx \reel^d\to[0,1]$ be a smooth function such that $h=1$ on
$V_1$ and $h=0$ on $V_2^c$. The coefficient $b_\Phi(y)$ of \eqref
{aphi} is defined on $V_3$. On the other hand, from
Assumption \ref{ascoeff}, there exists $\eps>0$ such that
\[
|\lambda|\le\eps\quad\Rightarrow\quad a(U_2,\lambda)\subset U_3.
\]
Then $a_\Phi(y,\lambda)$ is well defined for $y\in V_2$ and $|\lambda
|\le\eps$, and takes its values in $V_3$. Thus,
for $y\in\reel^d$ and $|\lambda|\le\eps$, we can define
%
%e13 #&#
\begin{equation}\label{atphi}
 \bigl(\atilde_\Phi(y,\lambda),\btilde_\Phi(y)
\bigr)= %
\cases{ \bigl(h(y)a_\Phi(y,\lambda)+\bigl(1-h(y)
\bigr)y,h(y)b_\Phi(y)\bigr) &\quad if $y\in V_2$,
\cr
(y,0)&\quad otherwise}
\end{equation}
which is an interpolation between $(a_\Phi,b_\Phi)$ on $V_1$ and the
motionless process on the complement of $V_2$.
Notice that from Assumption \ref{ascoeff},
\[
\atilde_\Phi(y,\lambda)=y+h(y)\abar_\Phi(y)\lambda+\mathrm{O}
\bigl(|\lambda |^\gamma\bigr)
\]
for $\abar_\Phi$ given by \eqref{aphi}, and similarly for its
derivatives. We can consider on $\reel^d$ the equation
\eqref{eqeuclid} with coefficients $(\atilde_\Phi,\btilde_\Phi)$,
driven by
%
%e14 #&#
\begin{equation}\label{lameps}
\Lambda_t^\eps=\Lambda_t-\sum
_{s\le t}\Delta\Lambda_s1_{\{|\Delta
\Lambda_s|>\eps\}}.
\end{equation}
Our smoothness assumptions on $a$ and $b$ imply that it has a unique
solution $Y_t$ for $Y_0=\Phi(x_0)$ fixed. Let
\[
\tau=\inf\bigl\{t\ge0;Y_t\notin V_1\mbox{ or }|\Delta
\Lambda_t|>\eps\bigr\}.
\]
Defining $X_t=\Phi^{-1}(Y_t)$ for $t<\tau$, and $X_\tau=a(X_{\tau
-},\Delta\Lambda_\tau)$, the process $X$ is solution
of our equation \eqref{eqx} up to the time $\tau$; if $|\Delta
\Lambda_\tau|>\eps$ and $X_\tau\in U_1$, then we can
solve again the equation from this time $\tau$ and point $X_\tau$;
since jumps greater than $\eps$ are in finite number
on any finite time interval, we deduce the existence of a solution
$X_t$ up to its first exit time from $U_1$.
Conversely, we can go from $X$ to $Y$ and deduce the uniqueness of $X$
from the uniqueness of $Y$. Thus, equation~\eqref{eqx} has a unique solution up to the first exit time from $U_1$.

Let us deduce the existence of a solution $(X_t;0\le t\le t_0)$ for a
fixed $t_0>0$. We consider a locally finite
atlas of $M$ made of subsets $U_0(k)$ (any compact subset intersects
finitely many of these $U_0(k)$) such that
$U_0(k)$ is relatively compact in an open subset $U(k)$ of the type of
the set $U_1$ of the first part of the proof. We
can choose for any $x$ an index $k(x)$ such that $k$ is measurable and
$x\in U_0(k(x))$, and we solve the equation
starting from $x$ up to the first exit time $\tau_1$ of the set
$U(k(x))$ (apply the first part of the proof). Let
$\proba_x$ be the law of this solution. If $\delta$ is some
Riemannian distance on $M$, one can check that $\sup_{t\le
u}\delta(x,X_t)$ under $\proba_x$ converges in probability to 0 as
$u\downarrow0$, uniformly for $x$ such that $k(x)=k$
(this follows from the similar property satisfied by $Y=\Phi(X)$);
thus there exists $u_k>0$ such that
%
%e15 #&#
\begin{equation}\label{pxtauk}
\proba_x[\tau_1\ge u_k]\ge
\proba_x\Bigl[\sup_{t\le u_k}\delta (x,X_t)<\delta
\bigl(U_0(k),U(k)^c\bigr)\Bigr]\ge1/2
\end{equation}
if $k=k(x)$.

The equation is then solved by means of the following iterative
procedure. For a fixed initial condition $x_0$, we
solve the equation from time $\tau_0=0$ up to the exit time $\tau_1$
from $U(k(x_0))$. If $\tau_1\ge t_0$ (the time up
to which we want to solve the equation), we have obtained the solution
up to time $t_0$ and we can stop the procedure;
if $X_{\tau_1}=\infty$ the process stays at $\infty$ and the
procedure can also be stopped; otherwise, starting at time
$\tau_1$ from $X_{\tau_1}$, we solve the equation up to the exit time
$\tau_2$ from $U(k(X_{\tau_1}))$, and so on. We
stop the procedure either when $\tau_j\ge t_0$, or when the process
has been killed (jump to $\infty$). Thus, the
procedure goes on forever when $\tau_j<t_0$ and $X_{\tau_j}\ne\infty
$ for any $j$. On the other hand, denoting by
$\tribuf_t$ the filtration of $\Lambda_t$, we have from \eqref
{pxtauk} and the strong Markov property that
\[
\proba[\tau_{j+1}-\tau_j\ge u_k|
\tribuf_{\tau_j}]\ge1/2\qquad \mbox{on } A_j^k=\{
\tau_j<t_0\}\cap\bigl\{k(X_{\tau_j})=k\bigr\}.
\]
We deduce that
\[
\sum_{j\ge0}\proba\bigl[A_j^k
\bigr]\le2\sum_{j\ge0}\proba\bigl[A_j^k
\cap\{\tau_{j+1}-\tau_j\ge u_k\}\bigr]
\le2(1+t_0/u_k)
\]
because there are at most $t_0/u_k$ disjoint intervals of length $\ge
u_k$ included in $[0,t_0]$. Thus, for $k$ fixed,
$A_j^k$ cannot be satisfied infinitely many times. We deduce that if
the procedure goes on forever, then
$k(X_{\tau_j})$ tends to infinity, so $X_{\tau-}=\infty$ for $\tau
=\lim\tau_j$ (the solution explodes at time~$\tau$).
In this case the solution is obtained by putting $X_t=\infty$ for
$t\ge\tau$.

The uniqueness can be proved by considering the first time $\tau$ at
which two solutions $X_t^1$ and $X_t^2$ differ,
and by supposing that $\tau<\infty$ with positive probability; then
$X_{\tau-}^1=X_{\tau-}^2$, so $X_\tau^1=X_\tau^2$
and the uniqueness of the solution in $U(k)$ for $k=k(X_\tau^1)$ leads
to a contradiction.
\end{pf}

%s3.2 #&#
\subsection{The case with only small jumps}

In Lemma \ref{exuniq}, we have worked out a construction of the
process by means of local charts. We now verify that
this construction also provides a smoothness result for the law on
these local charts. In all the proofs, we choose a
Riemannian metric on $M$, and the associated $C^\infty$ measure $\mathrm{d}x$
with respect to which we will consider the
densities of $M$-valued random variables. We shall study the solution
$X^\eps$ of equation \eqref{eqx} driven by the
L\'{e}vy process without its jumps greater than $\eps$ (the process
$\Lambda^\eps$ defined in \eqref{lameps}). The
parameter $\eps$ will be fixed and will be assumed to be small enough;
notice however that the constants involved in
the calculations will not be uniform in $\eps$; this does not cause
any difficulty because we shall never take the
limit as $\eps\downarrow0$.

%le2 #&#
\begin{lemma}\label{smflat}
Consider open subsets of $M$ with relatively compact inclusions
$U_0\Subset U_1\Subset U_2\Subset U_3$, such that $U_3$
is diffeomorphic to an open subset of $\reel^d$. Moreover, let $B$ be
a neighbourhood of the diagonal in $M\times M$. We
consider the process $X^\eps$ solution of equation \eqref{eqx} driven
by $\Lambda^\eps$. The surjectivity of $\abar(x)$
is only assumed on $U_2$. Let $\ell\ge0$. Then the following
properties hold true if $\eps$ is small enough.
\begin{enumerate}
\item The law of the process $X^\eps$ starting from $x_0\in U_1$ and
killed at the exit from $U_1$ has a $C^\ell$
density on $U_0$.
\item The density $x\mapsto q^\eps(t,x_0,x)$ of this killed process
and more generally its derivatives of order
$j\le\ell$ are uniformly dominated by $t^{-(d+j)/\alpha}$ for
$0<t\le t_0$, $x_0\in U_1$ and $x\in U_0$.
\item The density and its derivatives up to order $\ell$ are uniformly
bounded for $0<t\le t_0$ and
$(x_0,x)\in(U_1\times U_0)\setminus B$.
\end{enumerate}
\end{lemma}

\begin{pf}
We apply the construction given in the proof of Lemma \ref{exuniq},
denote $V_i=\Phi(U_i)$, and obtain the process
$Y_t$ solution of \eqref{eqeuclid} driven by $\Lambda^\eps$, with
coefficients $(\atilde_\Phi,\btilde_\Phi)$ given by
\eqref{atphi}. Then $X^\eps$ can be written as $X_t^\eps=\Phi^{-1}(Y_t)$ strictly before the first exit from $U_1$, so
the killed processes $X^\eps$ and $\Phi^{-1}(Y)$ coincide. The
nonkilled process $Y_t$ has of course not a smooth law
for any initial condition, since it is motionless out of $V_2$. It is
however possible to modify $Y$ by adding extra
independent noise in its equation, without modifying the killed
process. For $\lambda\in\reel^m$, $|\lambda|\le\eps$
and $\lambda'\in\reel^d$, we can replace $\atilde_\Phi$ by
\[
\atilde_\Phi\bigl(y,\lambda,\lambda'\bigr)= %
\cases{ h(y)a_\Phi(y,\lambda)+\bigl(1-h(y)\bigr) \bigl(y+
\lambda'\bigr)&\quad  if $y\in V_2$,
\cr
y+\lambda'& \quad
otherwise,}
\]
so that
\[
\atilde_\Phi\bigl(y,\lambda,\lambda'\bigr)=y+h(y)
\abar_\Phi(y)\lambda +\bigl(1-h(y)\bigr)\lambda'+\mathrm{O}\bigl(|
\lambda|^\gamma\bigr).
\]
Then, letting $\Lambda'$ be a $d$-dimensional L\'{e}vy process
independent of $\Lambda$, satisfying Assumption \ref{aslevy}
and with jumps bounded by $\eps$, we can solve the equation \eqref
{eqeuclid} with coefficients
$(\atilde_\Phi,\btilde_\Phi)$ and driven by $(\Lambda^\eps
,\Lambda')$. The advantage is that now the differential of
$\atilde_\Phi$ with respect to $(\lambda,\lambda')$ at $(0,0)$ is
now surjective onto $\reel^d$, uniformly in $y$;
moreover, if $\eps$ has been chosen small enough and if $|\lambda|$
and $|\lambda'|$ are $\le\eps$, the map
$y\mapsto\atilde_\Phi(y,\lambda,\lambda')$ is a diffeomorphism
from $\reel^d$ onto itself, and its Jacobian determinant
is uniformly positive. Consequently, we can apply Theorem 1 of \cite
{picard97}, and deduce that $Y_t$ has a smooth
density $y\mapsto p_Y(t,y_0,y)$, with derivatives satisfying
%
%e16 #&#
\begin{equation}\label{supdjp}
 \sup_{y_0,y}\bigl|D^jp_Y(t,y_0,y)\bigr|
\le C_j t^{-(d+j)/\alpha}.
\end{equation}
Moreover, it is proved in Lemma 2 of \cite{picsav00} that $p_Y$ and
its derivatives up to order $\ell$ are actually
bounded as $t\downarrow0$ if the number of jumps necessary to go from
$y_0$ to $y$ is large enough; thus, for $c>0$
fixed and if $\eps$ has been chosen small enough,
%
%e17 #&#
\begin{equation}\label{supunif}
 \sup\bigl\{\bigl|D^jp_Y(t,y_0,y)\bigr|;0<t
\le t_0,|y-y_0|\ge c\bigr\}\le C_j.
\end{equation}
If $\tau$ is the first exit time of $Y$ from $V_1$, we have from
\[
\esp_{y_0}\bigl[f(Y_t)1_{\{t<\tau\}}\bigr]=
\esp_{y_0}\bigl[f(Y_t)\bigr] -\esp_{y_0}
\bigl[f(Y_t)1_{\{t\ge\tau\}}\bigr]
\]
and the strong Markov property that the process $Y$ killed at $\tau$
has density
\[
q_Y(t,y_0,y)=p_Y(t,y_0,y)-
\esp_{y_0} \bigl[p_Y(t-\tau,Y_\tau,y)1_{\{t\ge
\tau\}}
\bigr].
\]
The first term is estimated from \eqref{supdjp} and \eqref{supunif},
and for the second one, we notice that
$Y_\tau\notin V_1$, and that the number of jumps necessary to go from
$V_1^c$ into $V_0$ is large if $\eps$ is small
enough, so $y\mapsto p_Y(t-\tau,Y_\tau,y)$ and its derivatives up to
order $\ell$ are bounded on $V_0$. We deduce the
smoothness of $q_Y$, and, by applying $\Phi^{-1}$, the smoothness of
the law of the process $X^\eps$ killed as well as
the estimates claimed in the lemma.
\end{pf}

%le3 #&#
\begin{lemma}\label{smjump}
Let $U$ be a relatively compact open subset of $M$. The surjectivity of
$\abar(x)$ is only assumed on the closure of
$U$. Consider again the solution $X^\eps$ of \eqref{eqx} driven by
$\Lambda^\eps$. If $\eps$ is small enough, then
$X_t^\eps$ has a $C^\ell$ density $p^\eps(t,x_0,x)$ on $U$ for any
$x_0\in M$; the density and more generally its
derivatives of order $j\le\ell$ are uniformly dominated by
$t^{-(d+j)/\alpha}$, for $0<t\le t_0$, $x_0$ in $M$ and $x$
in $U$.
\end{lemma}

\begin{pf}
It is sufficient to prove the result for $U=U_0$, for open subsets
$U_0\Subset U_1\Subset U_2\Subset U_3\Subset
U_4\Subset M$ such that $U_4$ is diffeomorphic to an open subset of
$\reel^d$, and $\abar(x)$ is surjective on $U_4$;
the subset $U$ of the lemma can indeed be covered by a finite number of
such sets $U_0$. Put $\tau_0=0$, and
\[
\tau_k'=\inf\bigl\{t\ge\tau_k;X_t^\eps
\notin U_3 \bigr\},\qquad \tau_{k+1}=\inf\bigl\{t\ge
\tau_k';X_t^\eps\in
U_2 \bigr\}.
\]
We can associate to this sequence of stopping times an expansion for
the law of $X_t^\eps$ on~$U_0$,
\begin{eqnarray*}
\proba_{x_0}\bigl[X_t^\eps\in \mathrm{d}x\bigr]&=&\sum
_{k=0}^\infty\proba_{x_0}
\bigl[X_t^\eps\in \mathrm{d}x,\tau_k\le t<
\tau_k'\bigr]
\\
&=&\sum_{k=0}^\infty\esp_{x_0}
\bigl[\proba_{x_0} \bigl[X_t^\eps\in \mathrm{d}x,t<
\tau_k'|\tribuf_{\tau_k} \bigr]1_{\{t\ge\tau_k\}}
\bigr]
\\
&=&\sum_{k=0}^\infty\esp_{x_0}
\bigl[Q^\eps\bigl(t-\tau_k,X_{\tau_k}^\eps
,\mathrm{d}x\bigr)1_{\{t\ge\tau_k\}} \bigr],
\end{eqnarray*}
where $\tribuf_t$ is the filtration of $\Lambda$ and $Q^\eps$ is the
transition kernel of the process $X^\eps$ killed
at the exit from $U_3$. From Lemma \ref{smflat}, this kernel has for
any $\ell$ a $C^\ell$ density $q^\eps$ on $U_0$ if
$\eps$ is small enough, so the law of $X_t^\eps$ is absolutely
continuous on $U_0$ with density
%
%e18 #&#
\begin{eqnarray}\label{pqtx}
 p^\eps(t,x_0,x)&=&\sum
_{k=0}^\infty\esp_{x_0} \bigl[q^\eps
\bigl(t-\tau_k,X_{\tau_k}^\eps,x\bigr)1_{\{t\ge\tau_k\}}
\bigr]
\nonumber
\\[-8pt]
\\[-8pt]
&=&q^\eps(t,x_0,x)+\sum_{k=1}^\infty
\esp_{x_0} \bigl[q^\eps\bigl(t-\tau_k,X_{\tau_k}^\eps,x
\bigr)1_{\{t\ge\tau_k\}} \bigr].\nonumber
\end{eqnarray}
We already know from Lemma \ref{smflat} that the first term (which is
0 if $x_0\notin U_3$) and its derivatives are
dominated by $t^{-(d+j)/\alpha}$. Moreover, we can choose $\eps$
small enough so that the process cannot jump from
$U_2^c$ into $U_1$, and in this case $X_{\tau_k}^\eps\notin U_1$ for
$k\ge1$; in particular $(X_{\tau_k}^\eps,x)$
remains out of a neighbourhood of the diagonal of $M\times M$ for $x\in
U_0$. Thus, $q^\eps(t-\tau_k,X_{\tau_k}^\eps,x)$
and its derivatives are uniformly bounded on $U_0$ for $k\ge1$ (third
assertion of Lemma \ref{smflat}). Thus, the proof
of the lemma is complete from
%
%e19 #&#
\begin{equation}\label{djptx}
 \bigl|D^jp^\eps(t,x_0,x) \bigr|\le
C_j \Biggl(1_{U_3}(x_0)t^{-(d+j)/\alpha} +\sum
_{k=1}^\infty\proba_{x_0} [t\ge
\tau_k ] \Biggr)
\end{equation}
as soon as we prove that the series converges and is bounded. We have
similarly to \eqref{pxtauk} that
$\proba_x[\tau_0'\ge u]\ge1/2$ for $x\in U_2$ if $u$ is small
enough, so by applying the strong Markov property of~$X$,
\[
\proba_{x_0} \bigl[\tau_k'>u|
\tribuf_{\tau_k} \bigr] \ge\proba_{x_0} \bigl[
\tau_k'-\tau_k\ge u|\tribuf_{\tau_k}
\bigr]\ge1/2
\]
for $k\ge1$ on $\{\tau_k<\infty\}$. Thus,
\begin{eqnarray*}
\proba_{x_0}\bigl[X_t^\eps\in U_3
\bigr]&\ge&\sum_{k=1}^\infty
\proba_{x_0} \bigl[\tau_k\le t<\tau_k'
\bigr] =\sum_{k=1}^\infty\esp_{x_0}
\bigl[\proba_{x_0} \bigl[\tau_k'>t|
\tribuf_{\tau_k} \bigr]1_{\{t\ge\tau_k\}} \bigr]
\\
&\ge&\frac12\sum_{k=1}^\infty
\proba_{x_0}[t\ge\tau_k]
\end{eqnarray*}
for $t\le u$, so the series in \eqref{djptx} is bounded, and we have
\[
\bigl|D^jp^\eps(t,x_0,x) \bigr|\le C_j
\bigl(1_{U_3}(x_0)t^{-(d+j)/\alpha}+2\proba_{x_0}
\bigl[X_t^\eps\in U_3\bigr] \bigr)
\]
for $t\le u$. If $t>u$, we use the Markov property, write
\[
p^\eps(t,x_0,x)=\esp_{x_0}
\bigl[p^\eps\bigl(u,X_{t-u}^\eps,x\bigr) \bigr]
\]
and deduce
%
%e20 #&#
\begin{eqnarray}\label{estdjp}
\bigl|D^jp^\eps(t,x_0,x) \bigr|&\le& C_j
\bigl(\proba_{x_0}\bigl[X_{t-u}^\eps\in
U_3\bigr]t^{-(d+j)/\alpha}+2 \proba_{x_0}
\bigl[X_t^\eps\in U_3\bigr] \bigr)
\nonumber
\\[-8pt]
\\[-8pt]
&\le& C_j' t^{-(d+j)/\alpha}\proba_{x_0}
\bigl[X_t^\eps\in U_4\bigr]\nonumber
\end{eqnarray}
by using the fact that $\proba[X_t^\eps\in U_4|X_{t-u}^\eps]\ge1/2$
on $\{X_{t-u}^\eps\in U_3\}$ if $u$ has been chosen
small enough.
\end{pf}

We have proved the estimate \eqref{estdjp} which is more precise than
the statement of the lemma. This property will be
used in Section \ref{liehom}.

The absolute continuity of the law of $X_t$ claimed in Theorem \ref
{mainth} follows easily from Lemma \ref{smjump}, by
conditioning on the last jump of $\Lambda$ before time $t$ greater
than $\eps$. If $\tau$ is this last jump (we put
$\tau=0$ when there is no big jump), then we deduce that $X$ has a
density given by
%
%e21 #&#
\begin{equation}\label{abscont}
 p(t,x_0,x)=\esp_{x_0}
\bigl[p^\eps(t-\tau,X_\tau,x) \bigr].
\end{equation}
However, this formula is not sufficient to obtain the smoothness and
even the local boundedness of $p$, because
$p^\eps(t-\tau,X_\tau,x)$ is of order $(t-\tau)^{-d/\alpha}$, at
least when $X_\tau$ and $x$ are close to each other,
and $(t-\tau)^{-d/\alpha}$ is not integrable if $d\ge\alpha$.

%s3.3 #&#
\subsection{The case with big jumps}

In Lemma \ref{smjump}, we have proved the smoothness of the law when
big jumps of $\Lambda$ have been removed. We now
have to take into account the effect of these big jumps. Notice that
the following lemma completes the proof of Theorem
\ref{mainth} when $M$ is compact (choose $U=M$).

%le4 #&#
\begin{lemma}\label{bjxy}
Let $U$ be a relatively compact open subset of $M$. There exists a
Markov process $Y_t$ such that the laws of $X$ and
$Y$ killed at the exit from $U$ coincide, and $Y_t$ has on $U$ a
$C^\ell$ density $p_Y(t,y_0,y)$, the derivatives of
which satisfy
%
%e22 #&#
\begin{equation}\label{djpy}
\bigl|D^jp_Y(t,y_0,y)\bigr |\le
C_j t^{-(d+j)/\alpha}
\end{equation}
for $0<t\le t_0$, $y_0$ in $M$ and $y$ in $U$.
\end{lemma}

\begin{pf}
For any $\eps>0$, we have a decomposition $\gener=(\gener-\gener{\,^\eps})+\gener{\,^\eps}$ of the infinitesimal generator,
where $\gener-\gener^{\,\eps}$ is the generator of the process $X^\eps$
driven by $\Lambda^\eps$, and
\[
\gener^{\,\eps} f(x)=\int_{\{|\lambda|>\eps\}} \bigl(f\bigl(a(x,
\lambda)\bigr)-f(x) \bigr)\mu(\mathrm{d}\lambda).
\]
Choose $U\Subset U'\Subset U''\Subset M$. Fix a Riemannian metric on
$M$ and consider the Riemannian exponential
function; then its inverse\vadjust{\goodbreak} $\exp_x^{-1}y\in T_xM$ is well defined and
smooth if $y$ is close to $x$. Thus, if $\eps$ is
small enough, we can consider $a_0(x,\lambda)=\exp_x^{-1}a(x,\lambda
)$ for $x\in U''$ and $|\lambda|\le\eps$. Let
$h_0\dvtx M\to[0,1]$ be a smooth function such that $h_0=1$ on $U'$ and
$h_0=0$ on ${U''}^c$, and let
\[
\atilde(x,\lambda)=\exp_x \bigl(h_0(x)a_0(x,
\lambda) \bigr),\qquad \btilde(x)=h_0(x)b(x)
\]
for $|\lambda|\le\eps$. Let $\genertilde^{\;\eps}$ be the
infinitesimal generator $\gener-\gener^{\,\eps}$ where coefficients
$(a,b)$ have been replaced by $(\atilde,\btilde)$; this corresponds
to a process $\Xtilde$ driven by $\Lambda^\eps$,
which is interpolated between $X^\eps$ and the motionless process
(this is similar to the construction of Lemma
\ref{exuniq} but we are here on the manifold instead of $\reel^d$);
it satisfies Assumption \ref{ascoeff}, except the
surjectivity condition which does not hold out of $U''$ but holds on
the closure of $U'$. On the other hand, let
$h\dvtx M\to[0,1]$ be a smooth function such that $h=1$ on $U$ and $h=0$ on
${U'}^c$, and define
\[
\generbar^{\;\eps} f(x)=\int_{\{|\lambda|>\eps\}} \bigl(f\bigl(a(x,
\lambda )\bigr)h\bigl(a(x,\lambda)\bigr)h(x) -f(x) \bigr)\mu(\mathrm{d}\lambda).
\]
This means that we consider the jumps $\lambda$ of $\Lambda$ greater
than $\eps$; if the process is at a point $x$
before this jump, we kill it with probability $1-h(x)$; if it is not
killed, it jumps to $x_1=a(x,\lambda)$ and is
killed with probability $1-h(x_1)$. We let $Y$ be the process with
generator $\genertilde^{\;\eps}+\generbar^{\;\eps}$; this is
the process $\Xtilde$ interlaced with jumps described by $\generbar^{\;\eps}$. This process enters our framework from
Section \ref{jcoeff}, and $X$ and $Y$ coincide when killed at the exit
from $U$.

From Lemma \ref{smjump}, if $\eps$ is small enough, the process
$\Xtilde$ has on $U'$ a $C^\ell$ density
$\ptilde(t,x_0,x)$ with respect to the Riemannian measure. On the
other hand, there are $N_t$ jumps
$\Delta\Lambda_{\tau_k}$ greater than $\eps$ on the time interval
$[0,t]$, for a random $\tau=(\tau_1,\tau_2,\ldots)$;
we also let $\tau_0=0$ and append a last $\tau_{N_t+1}=t$. Let $K$ be
the random index $k$ such that
$\tau_{k+1}-\tau_k$ is maximal. Then
\[
\proba[Y_{\tau_{K+1}-}\in \mathrm{d}y|\tau;\Lambda_s,0\le s\le
\tau_K ]=\ptilde(\tau_{K+1}-\tau_K,Y_{\tau_K},y)\, \mathrm{d}y
\]
on $U'$, with
%
%e23 #&#
\begin{equation}\label{djptau}
 \bigl|D^j\ptilde(\tau_{K+1}-\tau_K,Y_{\tau_K},y)
\bigr|\le C_j (\tau_{K+1}-\tau_K )^{-(d+j)/\alpha}.
\end{equation}
On $\{K=N_t\}$, we obtain the conditional density of $Y_t$. Otherwise,
we have to apply the jump at $\tau_{K+1}$ to
this distribution; we first kill the process with probability $1-h(y)$
and therefore get a $C^\ell$ law on $M$
supported by $U'$; from Assumption \ref{ascoeff}, this law is then
transformed by $a_\star^\eps$ into a $C^\ell$ law on
$M$, which is restricted into a $C^\ell$ law supported by $U'$ by the
second killing. We therefore obtain
\[
\proba[Y_{\tau_{K+1}}\in \mathrm{d}y|\tau;\Lambda_s,0\le s\le
\tau_K ]=p_\star(y)\,\mathrm{d}y
\]
for a conditional density $p_\star$ which is $C^\ell$, with
derivatives dominated as in \eqref{djptau}.

This density is then propagated from $\tau_{K+1}$ to $(\tau_{K+2})-$
by means of the semigroup of $\Xtilde$ with
generator $\genertilde^{\;\eps}$; if $\eps$ is small enough, then
$x\mapsto\atilde(x,\lambda)$ are diffeomorphisms of $M$
onto itself for $|\lambda|\le\eps$, and the process $\Xtilde$ can
be written as $\Xtilde_t=\Phi_t(\Xtilde_0)$ for a
flow of diffeomorphisms $\Phi_t$ of $M$ onto itself; the technique of
\cite{fujiwara91} for compact manifolds can be
adapted to our case since the process is motionless out of a compact
part of $M$. We choose a copy of $\Phi$ which is
independent of $\Lambda$, and obtain
%
%e24 #&#
\begin{eqnarray}\label{prytau}
&&\proba[Y_{\tau_{K+2}-}\in \mathrm{d}y|\tau;\Lambda_s,0\le s
\le\tau_K ]
\\
&&\quad=\esptilde\bigl[p_\star\bigl(\Phi_{\tau_{K+2}-\tau_{K+1}}^{-1}(y)
\bigr)\bigl |\det \bigl(D\Phi_{\tau_{K+2}-\tau_{K+1}}^{-1}\bigr) (y) \bigr| \bigr]\,\mathrm{d}y,
\nonumber
\end{eqnarray}
where $\esptilde$ is the expectation only with respect to $\Phi$, and
the determinant is computed relatively to
orthonormal bases on the tangent spaces. The differential of $\Phi_t$
is solution of
\[
D\Phi_{t+\mathrm{d}t}(x)=D\atilde\bigl(\Phi_t(x),\mathrm{d}
\Lambdatilde_t\bigr)D\Phi_t(x)+D\btilde\bigl(
\Phi_t(x)\bigr)D\Phi_t(x)\,\mathrm{d}t, \qquad D\Phi_0(x)=I.
\]
If $M$ is embedded in $\reel^N$, this can be transformed into an
equation on $\reel^N$ and it is a standard procedure
to prove that $\sup_{t\le T}|D\Phi_t(x)|$ has bounded moments,
uniformly in $x$: use the technique of
\cite{fujikunita85}. By differentiating this equation, the same
property holds true for higher order derivatives, so
actually $\sup_{t,x}|D^j\Phi_t(x)|$, $t\le t_0$, has bounded moments.
The same property can be verified for the
derivatives of the inverse map $\Phi_t^{-1}$, by looking at the
equation of its derivative. Thus, in \eqref{prytau}, we
obtain an estimate on the derivatives of the conditional density of
$Y_{\tau_{K+2}-}$ similar to \eqref{djptau}. By
iterating this procedure on all the subsequent jumps $\tau_k$, we
prove that $Y_t$ has a conditional density; the
unconditioned density $p_Y$ is then obtained by taking the expectation
and satisfy
\[
\bigl|D^jp_Y(y) \bigr|\le C \esp\bigl[\mathrm{e}^{CN_t}(
\tau_{K+1}-\tau_K)^{-(d+j)/\alpha} \bigr]
\]
on $U'$, because each jump at $\tau_k$ and each use of the flow $\Phi
$ between $\tau_k$ and $\tau_{k+1}$ appends a
multiplicative constant in the estimation, and they are at most $N_t$
of these jumps. The number $N_t$ of big jumps is
a Poisson variable so has finite exponential moments. Moreover,
\[
\tau_{K+1}-\tau_K\ge\frac{t}{N_t+1},
\]
so
\[
(\tau_{K+1}-\tau_K)^{-(d+j)/\alpha}\le(N_t+1)^{(d+j)/\alpha
}t^{-(d+j)/\alpha}
\]
and we can conclude.
\end{pf}

We now give an estimation of the probability for $X_t$ to be in some
relatively compact open subset when one needs many
jumps to come from the initial condition.

%le5 #&#
\begin{lemma}\label{unvn}
Consider open subsets of $M$ with relatively compact inclusions
$U_n\Subset V_n\Subset U_{n+1}$, and suppose that $X$
cannot jump from $U_{n+1}^c$ into $V_n$. Let $n\ge1$. Then $\proba_{x_0}[X_t\in U_0]$ is $\mathrm{O}(t^n)$ as $t\downarrow0$
uniformly for $x_0$ in $V_{n-1}^c$.
\end{lemma}

\begin{pf}
A similar result was proved in Lemma 1 of \cite{picsav00} for $\reel^d$; the proof of this variant is much simpler.
Let $h_k$, $0\le k\le n-1$, be a smooth function with values in
$[0,1]$, which is 1 on $U_k$ and 0 on $V_k^c$. Let
$C_k=\sup|\gener h_k|$ for the generator $\gener$ of $X$. Then
$h_k(x_0)=0$ for $x_0\in V_{n-1}^c\subset V_k^c$, and
$\gener h_k(x)=0$ for $x\in U_{k+1}^c$, so
\[
\proba_{x_0}[X_t\in U_k]\le
\esp_{x_0}\bigl[h_k(X_t)\bigr]=
\esp_{x_0}\int_0^t\gener
h_k(X_s)\,\mathrm{d}s \le C_k t\sup_{0\le s\le t}
\proba_{x_0}[X_s\in U_{k+1}].
\]
Applying this inequality for $k=0,1,\ldots,n-1$ completes the proof.
\end{pf}

\begin{pf*}{Proof of Theorem \ref{mainth}}
The absolute continuity has already been proved in \eqref{abscont},
the smoothness has been obtained in Lemma
\ref{bjxy} in the compact case, so we now have to study the smoothness
of the density in the non compact case. Let
$U_0$ be a relatively compact open subset of $M$ which is diffeomorphic
to an open subset of $\reel^d$, and let
$K\subset U_0$ be compact. Under Assumption \ref{asjump}, there exists
a sequence $(U_n,V_n)$ satisfying the conditions
of Lemma \ref{unvn} starting from the given $U_0$. For $\eps$ and $n$
which will be chosen respectively small and large
enough, there also exists a process $Y$ constructed in Lemma \ref
{bjxy} for $U=U_{n+1}$, with a $C^\ell$ density
$p_Y(t,y_0,y)$ on $U_{n+1}$. The law of $Y$ killed at the exit $\tau$
from $U_{n+1}$ has on $U_0$ the density
%
%e25 #&#
\begin{equation}\label{qypy}
 q_Y(t,y_0,y)=p_Y(t,y_0,y)-
\esp_{y_0} \bigl[p_Y(t-\tau,Y_\tau,y)1_{\{\tau
<t\}}
\bigr].
\end{equation}
Consider the process $Y$ when it is out of $U_{n+1}$; its big jumps
(coming from $\generbar^{\;\eps}$) are included in
jumps of $X$, so by construction of the sequence $(U_n,V_n)$, it is not
possible to jump directly from $U_{n+1}^c$ into
$V_n$; the small jumps (coming from $\genertilde^{\;\eps}$) have been
modified, but they are small, so if $\eps$ is small
enough, it is again not possible to jump directly from $U_{n+1}^c$ into
$V_n$. On the other hand, $Y$ coincides with
$X$ on $U_{n+1}$ and has therefore the same jumps, so it cannot jump
from $U_{k+1}^c$ into $V_k$ for $k<n$. Thus, we can
apply Lemma \ref{unvn} to $Y$ and deduce that $\proba_{y_0}[Y_t\in
U_0]$ is $\mathrm{O}(t^n)$ for $y_0\notin V_{n-1}$; we also
have the uniform estimates \eqref{djpy}, and we can deduce as in Lemma
2 of \cite{picsav00} that if $n$ has been
chosen large enough, then $D^jp_Y(t,y_0,y)$ is uniformly bounded for
$0<t\le t_0$, $y\in K$ and $y_0\notin V_{n-1}$
(this result was proved on $\reel^d$ but $U_0$ can be viewed as a
subset of $\reel^d$). These estimates on $p_Y$ imply
by means of \eqref{qypy} that $D^jq_Y(t,y_0,y)$ is $\mathrm{O}(t^{-(d+j)/\alpha
})$ uniformly for $0<t\le t_0$, $y_0\in U_{n+1}$,
$y\in K$, and is bounded if $y_0\in U_{n+1}\setminus V_{n-1}$.

The killed processes $X$ and $Y$ coincide, and the smoothness of the
non killed process $X$ is then deduced as in the
proof of Lemma \ref{smjump} by considering successive exits from
$U_{n+1}$ and entrances into $V_n$; we obtain an
expansion similar to \eqref{pqtx}. When the process enters $V_n$, it
is at a point out of $V_{n-1}$, so we have the
uniform boundedness of the derivatives of $q_Y$ starting from this
point, and we can proceed and estimate the series as
in Lemma~\ref{smjump}.
\end{pf*}

\begin{pf*}{Proof of Theorem \ref{smoothj}}
The decomposition $\mu=\mu^\flat+\mu^\sharp$ of the L\'{e}vy
measure corresponds to a decomposition
$\Lambda_t=\Lambda_t^\flat+\Lambda_t^\sharp$ of the L\'{e}vy
process into independent L\'{e}vy processes, where
$\Lambda^\sharp$ is of pure jump type. We can apply Theorem \ref
{mainth} and deduce that the process $X^\flat$ driven
by $\Lambda^\flat$ has a smooth density $p^\flat$. Let $\tau$ be
the time of the last jump of $\Lambda^\sharp$ before
$t$, with $\tau=0$ when $\Lambda^\sharp$ has no jump on $[0,t]$;
this last event has probability
$\exp-t\mu^\sharp(\reel^m)$. Similarly to \eqref{abscont}, the
density of $X_t$ can be written as
\begin{eqnarray*}
p(t,x_0,x)&=&\esp_{x_0} \bigl[p^\flat(t-
\tau,X_\tau,x) \bigr]
\\
&=&p^\flat(t,x_0,x)\exp\bigl(-t\mu^\sharp\bigl(
\reel^m\bigr) \bigr) +\esp_{x_0} \bigl[p^\flat(t-
\tau,X_\tau,x)1_{\{\tau>0\}} \bigr].
\end{eqnarray*}
The first term is smooth. On the event $\{\tau>0\}$, it follows from
the assumption of the theorem about $\mu_x^\sharp$
that the conditional law of $X_\tau$ given $\Lambda^\flat$ has a
smooth density $p_0$ on $M$, so
\[
\esp_{x_0} \bigl[p^\flat(t-\tau,X_\tau,x)1_{\{\tau>0\}}
\bigr] =\esp_{x_0} \biggl[1_{\{\tau>0\}}\int p^\flat(t-
\tau,z,x)p_0(z)\,\mathrm{d}z \biggr].
\]
We want to verify that the smoothness of $p_0$ is preserved by $p^\flat
(t-\tau,\cdot,\cdot)$. We know from the proof of Theorem
\ref{main} that $p^\flat(t-\tau,z,\cdot)$ is $C_b^\ell$ on $K$ if the
initial condition $z$ is out of a large enough
subset, so it is sufficient to consider the case of a smooth $p_0$ with
compact support. It also follows from previous
proof that it is sufficient to prove the result for the modified
process $Y$, and the propagation of the smoothness of
the law by the semigroup of $Y$ has been obtained in the proof of Lemma
\ref{bjxy}.
\end{pf*}

%s4 #&#
\section{Lie groups and homogeneous spaces}\label{liehom}

In Theorem \ref{mainth}, we have assumed that the process cannot come
from infinity by jumping (except in the case of
smooth jumps of Theorem \ref{smoothj}). This assumption is not needed
in the affine case $M=\reel^d$, but other
conditions are required in this case, and assumptions made for instance
in \cite{picard96} are not intrinsic if
$\reel^d$ is viewed as a differentiable manifold (they are not
invariant by diffeomorphisms). We can of course apply
the affine theorem when $M$ is diffeomorphic to $\reel^d$, but again
the conditions on the coefficients will depend on
the diffeomorphism.

Thus, if jumps are not bounded, we need additional structure on $M$ (as
this is the case on $\reel^d$ where the affine
structure is used). We consider here the case of L\'{e}vy processes on
Lie groups; more generally, we consider the case
where $M$ is an homogeneous space on which a Lie group $G$ acts, and
$X_t$ is the projection on $M$ of a L\'{e}vy process
on $G$.

Exposition about L\'{e}vy processes on Lie groups can be found in \cite
{liao04}. We shall use the theory of integration on
Lie groups or homogeneous spaces (Haar measures, invariant and
relatively invariant measures), which is explained in
\cite{bourbaki63}, see also for instance \cite{wijsman90}. We recall
here the points which are useful for our study.

%s4.1 #&#
\subsection{Lie groups}

Let $M=G$ be a $d$-dimensional Lie group with neutral element $e$ and
Lie algebra $\lieal$; as a vector space, $\lieal$
is the tangent space $T_eG$; it can be identified to the space of left
invariant vector fields on $G$; the Lie bracket
of two elements of $\lieal$, as well as the exponential map $\exp
\dvtx\lieal\to G$ can be constructed from this
identification. We can choose as a smooth reference measure on $G$ a
left Haar measure $\haar_\leftarrow^{\;G}$, or a right
Haar measure $\haar_\rightarrow^{\;G}$; each of them is unique modulo a
multiplicative constant; they satisfy
%
%e26 #&#
\begin{eqnarray}\label{haarlr}
\haar_\leftarrow^{\;G}(gA)&=&
\haar_\leftarrow^{\;G}(A),\qquad \haar_\leftarrow^{\;G}
\bigl(Ag^{-1}\bigr)=\chi_G(g)\haar_\leftarrow^{\;G}(A),\nonumber
\\[-8pt]
\\[-8pt]
\haar_\rightarrow^{\;G}(Ag)&=&\haar_\rightarrow^{\;G}(A),\qquad
\haar_\rightarrow^{\;G}(gA)=\chi_G(g)
\haar_\rightarrow^{\;G}(A)\nonumber
\end{eqnarray}
for a group homomorphism $\chi_G\dvtx G\to\reel_+^\star$ which is the
modulus of $G$. If we are given $\haar_\leftarrow^{\;G}$,
we can define $\haar_\rightarrow^{\;G}$ by
%
%e27 #&#
\begin{equation}\label{hgrl}
\haar_\rightarrow^{\;G}(\mathrm{d}g)=\chi_G(g)
\haar_\leftarrow^{\;G}(\mathrm{d}g).
\end{equation}
The group $G$ is said to be unimodular if $\chi_G\equiv1$; this holds
true when $G$ is compact because any group
homomorphism from $G$ into $\reel_+^\star$ must be equal to 1. Let
$\Ad\dvtx G\to GL(\lieal)$ be the adjoint representation
of $G$; this means that $\Ad_g$ is the differential at $e$ of the
inner automorphism $x\mapsto gxg^{-1}$. Then
%
%e28 #&#
\begin{equation}\label{chidet}
 \chi_G(g)= |\det\Ad_g |.
\end{equation}
On $G$ there are two differential calculi, a left invariant one and a
right invariant one. The left and right invariant
derivatives are the linear forms
\[
D_\leftarrow f(g)u=\frac{\mathrm{d}}{\mathrm{d}\eps}\biggl|_{\eps=0}f\bigl(g\exp(\eps u)
\bigr),\qquad D_\rightarrow f(g)u=\frac{\mathrm{d}}{\mathrm{d}\eps}\biggr|_{\eps=0}f\bigl(\exp(\eps
u)g\bigr)
\]
for smooth functions $f\dvtx G\to\reel$ and $u\in\lieal$. The invariance
means that $D_\leftarrow L_hf=L_h D_\leftarrow f$
and $D_\rightarrow R_hf=R_hD_\rightarrow f$, with the notations
$L_hf(g)=f(hg)$ and $R_hf(g)=f(gh)$. The left and right
invariant derivatives are related to each other by
%
%e29 #&#
\begin{equation}\label{dleftright}
D_\rightarrow f(g)u=D_\leftarrow f(g)
\Ad_g^{-1}u.
\end{equation}
If we choose an inner product on $\lieal$, we can consider the norms
$|D_\leftarrow f(g)|$ and $|D_\rightarrow f(g)|$
of the linear forms, and norms corresponding to different inner
products are of course equivalent. We can also consider
classes $C_{b,\leftarrow}^\ell$ or $C_{b,\rightarrow}^\ell$ of
functions for which the left or right invariant
derivatives are bounded (without boundedness we can simply use the
notation $C^\ell$ since the classes for the left and
right calculi coincide).

Let $X_t$ be a left L\'{e}vy process on $G$ with initial condition
$X_0=e$. This is a strong Markov process which is
invariant under left multiplication, so that its semigroup satisfies
$P_tL_h=L_hP_t$; equivalently, the infinitesimal
generator should satisfy $\gener L_h=L_h\gener$. For $0\le s\le t$,
the variable $X_s^{-1}X_t$ must be independent of
$(X_u;0\le u\le s)$ and must have the same law as $X_{t-s}$. We
consider here the subclass of L\'{e}vy processes without
Brownian part. Let $V$ be a relatively compact neighbourhood of $e$ in
$G$ which is diffeomorphic to a neighbourhood
$U$ of 0 in $\lieal$ by means of the Lie exponential function. Then a
left L\'{e}vy process without Brownian part is
characterised by a drift $\kappa\in\lieal$ and a L\'{e}vy measure
$\mu_X$ on $G\setminus\{e\}$ which integrates smooth
bounded functions $f$ such that $f(e)=D_\leftarrow f(e)=0$; the
infinitesimal generator of $X$ can be written in the
Hunt form as
%
%e30 #&#
\begin{eqnarray}\label{infglevy}
\gener f(g)&=&D_\leftarrow f(g)\kappa+\int
_V \bigl(f(gx)-f(g)-D_\leftarrow f(g)
\exp^{-1}x \bigr)\mu_X(\mathrm{d}x)
\\
&&{}+\int_{V^c} \bigl(f(gx)-f(g) \bigr)\mu_X(\mathrm{d}x).
\nonumber
\end{eqnarray}
It is explained in \cite{applebaumkuni93} that $X$ can be viewed as
the solution of an equation driven by a Poisson
measure on $\reel_+\times G$; by means of a technique similar to
Section \ref{jcoeff}, it is also the solution of an
equation of type \eqref{eqx} driven by a $\lieal$-valued L\'{e}vy
process. More precisely, let $\mu$ be the measure on
$U\setminus\{0\}$ which is the image of $\mu_X |_V$ by $\exp^{-1}$;
on the other hand, there exists a bi-measurable
bijection $i$ from $V^c$ onto a Borel subset $U'$ of $U^c$, so we can
let $\mu$ be on $U^c$ the image of
$\mu_X |_{V^c}$ by $i$. Then
\begin{eqnarray*}
\gener f(g)&=&D_\leftarrow f(g)\kappa+\int_U
\bigl(f(g\exp\lambda )-f(g)-D_\leftarrow f(g)\lambda\bigr)\mu(\mathrm{d}\lambda)
\\
&&{}+\int_{U'} \bigl(f\bigl(g i^{-1}(\lambda)
\bigr)-f(g) \bigr)\mu(\mathrm{d}\lambda).
\end{eqnarray*}
Thus $X_t$ can be viewed as the solution of $X_t=a(X_{t-},\mathrm{d}\Lambda_t)$, where $\Lambda_t$ is a L\'{e}vy process in the
vector space $\lieal$ with L\'{e}vy measure $\mu$ supported by $U\cup
U'$, and
\[
a(g,\lambda)= %
\cases{ g\exp\lambda & \quad if $\lambda\in U$,
\cr
g
i^{-1}(\lambda)& \quad if $\lambda\in U'$.}
\]
Then $\abar(g)\dvtx\lieal\to T_g(G)$ is the map which sends $u$ to the
value at $g$ of the left invariant vector field
associated to $u$, so it is bijective, and it is also not difficult to
verify the second part of Assumption
\ref{ascoeff}. The nondegeneracy condition of Assumption \ref{aslevy}
on $\mu$ is immediately transferred to an
assumption on $\mu_X$ as
%
%e31 #&#
\begin{equation}\label{nondlie}
 c\rho^{2-\alpha}|u|^2\le\int_{\{|\exp^{-1}x|\le\rho\}}
\bigl\langle \exp^{-1}x,u \bigr\rangle^2
\mu_X(\mathrm{d}x)\le C\rho^{2-\alpha}|u|^2
\end{equation}
for $u\in\lieal$ and $\rho\le\rho_0$ small enough so that $\exp^{-1}$ is well defined; the additional condition in the
case $\alpha=1$ is written as
%
%e32 #&#
\begin{equation}\label{symlie}
 \operatorname{\lim\sup}\limits_{\eps\downarrow0} \biggl|\int_{\{\eps\le|\exp^{-1}x|\le\rho
_0\}}
\exp^{-1}x\mu_X(\mathrm{d}x) \biggr|<\infty.
\end{equation}
These two conditions do not depend on the choice of the inner product
on $\lieal$. If $\alpha<1$, we also suppose that
the process is pure jump, so that \eqref{infglevy} becomes
\[
\gener f(g)=\int_G \bigl(f(gx)-f(g) \bigr)
\mu_X(\mathrm{d}x).
\]
For the first part of Assumption \ref{ascoeff}, we have to study the
propagation of the smoothness of the measure by a
right translation; this is easy if we choose a right Haar measure as a
reference measure.

Thus, under these conditions, we can apply Theorem \ref{mainth} and
deduce that $X_t$ has a smooth density if $\mu_X$
has compact support. If the support is not compact, a possibility is to
use Theorem~\ref{smjump}. Otherwise, we can use
the following result.

%th4 #&#
\begin{theorem}\label{llevy}
Let $X_t$ be a left L\'{e}vy process on $G$ with $X_0=e$, the L\'{e}vy
measure $\mu_X$ of which satisfies \eqref{nondlie}, and
the additional condition \eqref{symlie} if $\alpha=1$. If $\alpha
<1$, suppose moreover that $X_t$ is a pure jump
process. Then the law of $X_t$, $t>0$, is absolutely continuous with
respect to the left Haar measure
$\haar_\leftarrow^{\;G}$. Let $\ell\ge0$. If
%
%e33 #&#
\begin{equation}\label{chigad}
 \int_{V^c}\chi_G(g) |
\Ad_g |^j\mu_X(\mathrm{d}g)<\infty
\end{equation}
for a relatively compact neighbourhood $V$ of $e$ and for $j\le\ell$,
then the density is in $C_{b,\leftarrow}^\ell$.
In particular, if
%
%e34 #&#
\begin{equation}\label{adgj}
 \int_{V^c} |\Ad_g |^j
\mu_X(\mathrm{d}g)<\infty
\end{equation}
for any $j$, the density is in $C_{b,\leftarrow}^\infty$.
\end{theorem}

\begin{pf}
The absolute continuity follows from Theorem \ref{mainth}. Assume now
\eqref{chigad} for $j=0$. Let $V$ be a relatively
compact neighbourhood of $e$, and let $\mu_X^\flat$ and $\mu_X^\sharp$ be the restrictions of $\mu_X$ to $V$ and $V^c$;
as usually, $X$ can be written on $[0,t]$ as a L\'{e}vy process
$X^\flat$ with L\'{e}vy measure $\mu_X^\flat$ interlaced with
$N_t$ big jumps at times $\tau_k$ described by $\mu_X^\sharp$. We
proceed as in the proof of Lemma \ref{bjxy} and let
$\tau_{K+1}$ be the end of the longest subinterval of $[0,t]$ without
big jump. Then, from Theorem \ref{mainth},
conditionally on $(\tau_k)$, the variable $Y=X_{\tau_{K+1}-}$ has
with respect to the left Haar measure a $C^\ell$
density $p_Y$ satisfying
%
%e35 #&#
\begin{equation}\label{dleftpy}
\bigl |D_\leftarrow^jp_Y(y) \bigr|\le
C_j(\tau_{K+1}-\tau_K)^{-(d+j)/\alpha}
\end{equation}
on $V$, and the right-hand side is integrable as in Lemma \ref{bjxy}.
But Theorem \ref{mainth} also states that the
estimate is uniform with respect to the initial condition, so the same
estimate holds true for any $gY$, and one
deduces that \eqref{dleftpy} holds uniformly on $G$. Then, we have
$X_t=YZ$ where $Z=X_{\tau_{K+1}-}^{-1}X_t$ is,
conditionally on $(\tau_k)$, independent of $Y$. Conditionally on $Z$
and $(\tau_k)$, the variable $X_t$ is therefore
absolutely continuous with density
%
%e36 #&#
\begin{equation}\label{pxzchi}
 p(x|Z)=\chi_G(Z)p_Y\bigl(xZ^{-1}
\bigr)
\end{equation}
(this follows from \eqref{haarlr}). Thus, by integrating this formula
and applying \eqref{dleftpy}, we deduce that the
density of $X_t$ is bounded and continuous as soon as $\chi_G(Z)$ is
integrable. On the other hand,
\[
\chi_G(Z)=\chi_G \bigl(\bigl(X_{\tau_{K+1}}^\flat
\bigr)^{-1} X_t^\flat\bigr)\prod
_{k=K+1}^{N_t}\chi_G \bigl(X_{\tau_k-}^{-1}X_{\tau_k}
\bigr),
\]
where the different terms are conditionally independent given $(\tau_k)$; the process $\chi_G(X_t^\flat)$ is a
geometric L\'{e}vy process with bounded jumps, so the first term has
bounded conditional expectation, and we deduce from
\eqref{chigad} for $j=0$ that the conditional expectation of $\chi_G(Z)$ is bounded by some exponential of $N_t$, so
$\chi_G(Z)$ is integrable and the case $\ell=0$ is proved.

We can differentiate \eqref{pxzchi} and get
\[
D_\leftarrow p(x|Z)u=\chi_G(Z)D_\leftarrow
p_Y\bigl(xZ^{-1}\bigr)\Ad_Z(u)
\]
for $u\in\lieal$. For higher order derivatives, we have
\[
D_\leftarrow^jp(x|Z) (u_1,\ldots,u_j)=
\chi_G(Z)D_\leftarrow^jp_Y
\bigl(xZ^{-1}\bigr) \bigl(\Ad_Z(u_1),\ldots,
\Ad_Z(u_j) \bigr).
\]
We deduce the smoothness of the law of $X_t$ and the boundedness of its
derivatives if we prove that
$\chi_G(Z)|\Ad_Z|^j$ is integrable. To this end, we notice that
\begin{eqnarray*}
\chi_G(Z)|\Ad_Z|^j&\le&\prod
_{k=K+1}^{N_t} \bigl(\chi_G \bigl(
\bigl(X_{\tau
_k}^\flat\bigr)^{-1}X_{\tau_{k+1}}^\flat
\bigr) |\Ad_{ ((X_{\tau_k}^\flat)^{-1}X_{\tau_{k+1}}^\flat)} |^j \bigr)
\\
&&{}\times\prod_{k=K+1}^{N_t} \bigl(
\chi_G \bigl(X_{\tau_k-}^{-1}X_{\tau_k} \bigr)
| \Ad_{ (X_{\tau_k-}^{-1}X_{\tau_k} )} |^j \bigr)
\nonumber
\end{eqnarray*}
where the different terms are again conditionally independent given
$(\tau_k)$; the integrability is deduced similarly
to the case $j=0$ by using \eqref{chigad} for the terms of the second
product. The last claim of the theorem follows
from \eqref{chidet}.
\end{pf}

Condition \eqref{chigad} is trivially satisfied for $j=0$ when $G$ is
unimodular, so in this case we obtain the
existence of a continuous bounded density without any assumption on the
big jumps.

Notice also that \eqref{adgj} is related to the existence of
exponential moments for big jumps. Assume that $\mu_X$ is
the image of a measure $\mu$ on $\lieal$ by the Lie exponential (this
holds for instance when the exponential is
surjective). For $\lambda\in\lieal$, let ${\mathrm{ad}}_\lambda
\dvtx\lieal\to\lieal$ be the adjoint action given by
${\mathrm{ad}}_\lambda(u)=[\lambda,u]$ for the Lie bracket $[\cdot,\cdot]$.
We have $\Ad_{\exp(\lambda)}=\exp({\mathrm{ad}}_\lambda)$, so
\[
|\Ad_{\exp(\lambda)}|\le\exp\bigl(|{\mathrm{ad}}_\lambda|\bigr)\le\exp \bigl(c|
\lambda|\bigr),
\]
and \eqref{adgj} is satisfied if
\[
\int_{\{|\lambda|>1\}}\exp\bigl(C|\lambda| \bigr)\mu(\mathrm{d}\lambda)<\infty
\]
for any $C$. This is however the worst case. If $G$ is nilpotent, then
the expansion of $\exp({\mathrm{ad}}_\lambda)$ is finite,
and exponential moments can be replaced by ordinary moments; notice
however that the class of stable processes
introduced by \cite{kunita94} on simply connected nilpotent groups
does not enter our framework, because
\eqref{nondlie} is not satisfied.

Theorem \ref{llevy} can of course be translated into the case of right
Haar measure, right invariant derivatives, and
right L\'{e}vy processes (invariant by right multiplication). Then the
conditions \eqref{chigad} and \eqref{adgj} are
replaced by conditions on $\chi_G(g)^{-1}$ and $|\Ad_g^{-1}|$. On the
other hand, left and right L\'{e}vy processes with
the same infinitesimal generator at $e$ have the same law at any fixed
time $t$, because the right L\'{e}vy process $Y$ can
be deduced from the left process $X$ on $[0,t]$ by the formula
$Y_s=X_{t-s}^{-1}X_t$. Thus, in order to study the law
of $X_t=Y_t$, one can choose between left and right calculi. Notice
however that if for instance we apply the result
for the right L\'{e}vy process, we obtain that the density
$p_\rightarrow$ of $X_t$ with respect to $\haar_\rightarrow^{\;G}$
is of class $C_{b,\rightarrow}^\ell$; from \eqref{hgrl}, the
relation between the densities with respect to the left
and right Haar measures is $p_\rightarrow=\chi_G^{-1}p_\leftarrow$;
the density $p_\leftarrow$ is of class $C^\ell$,
but not necessarily of class $C_{b,\leftarrow}^\ell$ (the relation
between left and right invariant derivatives is
given in \eqref{dleftright}).

%s4.2 #&#
\subsection{Homogeneous spaces}\label{homog}

We now consider the case where the manifold $M$ is an homogeneous space
$M=G/H$. More precisely, $G$ is a
$m$-dimensional Lie group of transformations acting transitively and
smoothly on the left on $M$, and $H=\{h\in
G;h(o)=o\}$ is the isotropy group of some fixed point $o$ of $M$; the
projection $\pi\dvtx G\to M$ is given by
$\pi(g)=g(o)$. We can choose a Lebesgue measurable section $S$ of $\pi^{-1}$ (which exists from the measurable section
theorem), and any $g$ in $G$ can then be uniquely written as $g=S_xh$
for $x=\pi(g)=g(o)$ and some $h$ in $H$. The
action of $G$ on $M$ can be written as $g(y)=\pi(gS_y)$. We will
denote by $\lieal$ and $\lieh$ the Lie algebras of $G$
and $H$.

We can look for a measure on $M$ which would be invariant under the
action of $G$, but such a measure does not always
exist. We therefore weaken the invariance into a relative invariance
property, see \cite{bourbaki63,wijsman90} for an
introduction to the topic and some of the properties which are given
below; we say that a Radon non-identically zero
measure $\haar^{\;M}$ is relatively invariant under the action of $G$ with
multiplier $\chi$ if $\chi\dvtx G\to\reel_+^\star$ is
a group homomorphism and
\[
\haar^{\;M}\bigl(g(A)\bigr)=\chi(g)\haar^{\;M}(A).
\]
Then the measure is invariant if $\chi\equiv1$; this is necessarily
the case when $G$ is compact. For instance, we have
seen in \eqref{haarlr} that the right Haar measure on $G$ is
relatively invariant under the left multiplication with
multiplier $\chi_G$. A relatively invariant measure also does not
always exist, but it exists in more general
situations than invariant measures. It exists on $M=G/H$ if and only if
$\chi_G/\chi_H\dvtx H\to\reel_+^\star$ can be
extended to a group homomorphism $\chi\dvtx G\to\reel_+^\star$; in this
case there exist a relatively invariant measure with
multiplier $\chi$, and this measure is unique modulo a multiplicative
constant. In particular, an invariant measure
exists if and only if $\chi_G=\chi_H$ on $H$; this property holds in
particular when $H$ is compact.

The relationship between left Haar measures on $G$ and $H$ and
relatively invariant measures on $M$ is the following
one. If $\haar_\leftarrow^{\;H}$ is a left Haar measure on $H$ and if
$\chi\dvtx G\to\reel_+^\star$ is a group homomorphism,
then a Radon measure $\haar^{\;M}$ on $M$ is relatively invariant with
multiplier $\chi$ if and only if the measure
$\haar_\leftarrow^{\;G}$ defined on $G$ by
%
%e37 #&#
\begin{equation}\label{fghl}
\int_Gf(g)\haar_\leftarrow^{\;G}(\mathrm{d}g)=
\int\int_{M\times H}\frac{f}{\chi
}(S_xh)
\haar^{\;M}(\mathrm{d}x)\haar_\leftarrow^{\;H}(\mathrm{d}h)
\end{equation}
is a left Haar measure. For right Haar measures, \eqref{fghl} becomes
(use \eqref{hgrl})
%
%e38 #&#
\begin{equation}\label{fghr}
\int_Gf(g)\haar_\rightarrow^{\;G}(\mathrm{d}g)=
\int\int_{M\times H}f(S_xh)\frac
{\chi_G}\chi(S_x)
\haar^{\;M}(\mathrm{d}x)\haar_\rightarrow^{\;H}(\mathrm{d}h).
\end{equation}

If $\Xi$ is a $G$-valued variable with densities $p_\leftarrow$ and
$p_\rightarrow$ with respect to
$\haar_\leftarrow^{\;G}$ and $\haar_\rightarrow^{\;G}$, we deduce from
\eqref{fghl} and \eqref{fghr} that $\pi(\Xi)$ has
density
%
%e39 #&#
\begin{equation}\label{plr}
p(x)=\int_H\frac{p_\leftarrow}\chi(S_xh)
\haar_\leftarrow^{\;H}(\mathrm{d}h) =\frac{\chi_G}\chi(S_x)
\int_Hp_\rightarrow(S_xh)
\haar_\rightarrow^{\;H}(\mathrm{d}h)
\end{equation}
with respect to $\haar^{\;M}$. Estimates on $p$ can therefore be deduced
from estimates on $p_\leftarrow$ or
$p_\rightarrow$, but this is clearly simpler when $H$ is compact.

We now explain how left and right invariant differential calculi on $G$
can be transported to $M$. For the left
invariant calculus, notice that if $F$ is a smooth function on $M$, we
can differentiate $f=F\circ\pi$ which is defined
on $G$ and consider
%
%e40 #&#
\begin{equation}\label{dlfxu}
D_\leftarrow F(x)u=D_\leftarrow f(S_x)u
\end{equation}
for $u\in\lieal$; notice that the result is 0 if $u$ is in $\lieh$,
so the differential is actually a linear form on
the vector space $\lieal/\lieh$. The problem is that it depends on
the choice of the section $S$; if $S'$ is another
section and if we fix $x$, then $S_x'=S_xh$ for some $h\in H$, and
\[
D_\leftarrow f\bigl(S_x'\bigr)u=D_\leftarrow
f(S_x)\Ad_hu.
\]
If however $H$ is compact, then we can choose on $\lieal$ an $\Ad
(H)$-invariant inner product (this means that
$|\Ad_hu|=|u|$ for $h$ in $H$), we can consider $D_\leftarrow F(x)$ on
the orthogonal $\liep$ of $\lieh$ in $\lieal$,
and the norm $|D_\leftarrow F(x)|$ will not depend on the choice of
$S$. It is invariant under the action of $G$ in the
sense
\[
\bigl|D_\leftarrow(F\circ g) (x)\bigr |= \bigl|(D_\leftarrow F) \bigl(g(x)\bigr)\bigr |.
\]
Higher order derivatives have a similar behaviour, and we can consider
the classes of functions $C_{b,\leftarrow}^\ell$
on $M$.

For the right differential calculus, we only consider the case where
$(M,\cdot)$ is itself a Lie group, and left
translations of $M$ form a normal subgroup of $G$. In this case, a
canonical choice for the section $S$ is to let $S_x$
be the left translation by $x$, so that $S\dvtx M\to G$ is an injective
group homomorphism. The group $G$ is a semi-direct
product
%
%e41 #&#
\begin{equation}\label{gsmh}
 G=S(M)\rtimes H
\end{equation}
satisfying the commutation property $hS_x=S_{h(x)}h$ for $h\in H$. All
elements of $G$ can be written in the form
$S_xh$ for some $x\in M$ and $h\in H$, and we have the product rule
\[
S(x_1)h_1S(x_2)h_2=S
\bigl(x_1\cdot h_1(x_2)\bigr)h_1h_2.
\]
The vector space $\lieal$ can be written as $\lieal=\lieh\oplus
\liem$ for the Lie algebras $\lieh$ and $\liem$ of $H$
and $M\sim S(M)$. We can consider the differential $D_\rightarrow
F(x)u$ computed on the Lie group $M$ for $u\in\liem$,
and therefore the classes of functions $C_{b,\rightarrow}^\ell$. The
relation with the differential on $G$ is given by
\[
D_\rightarrow F(x)u=D_\rightarrow f(S_x)u\qquad \mbox{for }f=F
\circ\pi
\]
similar to \eqref{dlfxu}, because $\exp(\eps u)S_x=S_{\exp(\eps
u)x}$. The behaviour under the action of $G$ is
%
%e42 #&#
\begin{equation}\label{drfg}
 D_\rightarrow(F\circ g) (x)u=(D_\rightarrow F) \bigl(g(x)
\bigr)\Ad_gu.
\end{equation}
If $\haar_\rightarrow^{\;H}$ and $\haar_\rightarrow^{\;M}$ are right Haar
measures on $H$ and $M$, a right Haar measure can be
defined on $G$ by
%
%e43 #&#
\begin{equation}\label{fgrg}
\int_Gf(g)\haar_\rightarrow^{\;G}(\mathrm{d}g)=
\int\int_{M\times H}f(S_xh)\haar_\rightarrow^{\;M}(\mathrm{d}x)
\haar_\rightarrow^{\;H}(\mathrm{d}h)
\end{equation}
because the right-hand side is invariant under the right action of $H$
and $S(M)$ (use $S_xhS_y=S_{x\cdot h(y)}h$). By
comparing with \eqref{fghr}, we see that $\haar_\rightarrow^{\;M}$ is a
measure on $M$ with multiplier
%
%e44 #&#
\begin{equation}\label{multipm}
\chi(S_xh)=\chi_G(S_xh)/
\chi_H(h),
\end{equation}
since this is a group homomorphism and $\chi_G/\chi(S_x)=1$. In
particular, the Haar modulus of $M$ is
$\chi_M(x)=\chi(S_x)=\chi_G(S_x)$ and \eqref{plr} becomes
%
%e45 #&#
\begin{equation}\label{phpr}
p(x)=\int_Hp_\rightarrow(S_xh)
\haar_\rightarrow^{\;H}(\mathrm{d}h).
\end{equation}

In the next subsections, we study some classes of Markov processes on
$M$. In the first case, we assume that $H$ is
compact, and consider the class of Markov processes on $M$ which are
invariant under the action of $G$; they can be
written as the projection of left L\'{e}vy processes on $G$, and we
apply the left differential calculus. We also consider
the case where $X$ is the projection of a right L\'{e}vy process. In
the second case, we assume that $G$ is a semi-direct
product of type \eqref{gsmh}, and we let $X$ be the projection of a
right L\'{e}vy process on $G$; we then apply the right
differential calculus.

%s4.3 #&#
\subsection{Case 1: Compact isotropy subgroup}

We here assume that $M=G/H$ for a compact Lie subgroup $H$ of $G$; in
particular, there is on $H$ a unique probability
measure $\haar^{\;H}$ which is both a left and right Haar measure, and
there is on $M$ a measure $\haar^{\;M}$ which is
invariant under the action of $G$, and which is related to a left Haar
measure on $G$ by \eqref{fghl} with
$\chi\equiv1$. We choose an $\Ad(H)$-invariant inner product on
$\lieal$; it induces a $H$-invariant inner product on
the tangent space $T_oM$, and a $G$-invariant Riemannian metric on $M$.
We consider the class of Markov processes $X_t$
which are invariant under the action of $G$; this means that
$P_t(F\circ g)=(P_tF)\circ g$ for any $g\in G$; in
particular, the law of $(X_t;t\ge0)$ with initial condition $X_0=o$ is
invariant under the action of $H$. From Theorem
2.2 of \cite{liao04}, these processes are obtained as $X_t=\pi(\Xi_t)$, where $\Xi_t$ is a left L\'{e}vy process on $G$
which is invariant under the right action of $H$. The invariance of $X$
implies in particular that its L\'{e}vy kernel
$\mu_x$ can be deduced from $\mu_o$ which will be simply called the
L\'{e}vy measure $\mu_X$ of $X$; it is a $H$-invariant
measure on $M\setminus\{o\}$ which integrates $\delta^2(o,x)\wedge1$
(for the Riemannian distance $\delta$). The L\'{e}vy
process $\Xi$ on $G$ can be obtained by taking a L\'{e}vy measure
given by
%
%e46 #&#
\begin{equation}\label{muxia}
\mu_\Xi(A)=\int\int_{H\times M}1_A
\bigl(hS_xh^{-1}\bigr)\haar^{\;H}(\mathrm{d}h)
\mu_X(\mathrm{d}x),
\end{equation}
where the $H$-invariance of $\mu_X$ implies that $\mu_\Xi$ does not
depend on the choice of the section $S$ (see
\cite{liao04}). Then the generator of $\Xi$ can be written by means
of \eqref{infglevy}, for a neighbourhood $V_G=\exp
U_G$, where $U_G$ is a small enough ball of $\lieal$ (in particular
$U_G$ is $\Ad(H)$-invariant), and for a drift
$\kappa\in\lieal$ which is $\Ad(H)$-invariant.

For $0\le s\le t$ fixed, notice that $\Xi_t$ has the same law as $\Xi_s\Xi_{t-s}'$, for $\Xi'$ an independent copy of
$\Xi$; moreover, the variable $\Xi_s$ can be written as $S_{X_s}h_s$,
for $h_s$ a $H$-valued variable. Thus, letting
$X'=\pi(\Xi')$,
%
%e47 #&#
\begin{equation}\label{eqhom}
 X_t=\Xi_t(o)\sim\Xi_s
\Xi_{t-s}'(o)=S_{X_s}h_s
\bigl(X_{t-s}'\bigr)\sim S_{X_s}
\bigl(X_{t-s}'\bigr)
\end{equation}
because the law of $X_{t-s}'$ is $H$-invariant. More generally, if $\mu_1$ and $\mu_2$ are two $H$-invariant measures
on $M$, one can define the convolution $\mu_1\ast\mu_2$ as the image
of the product measure by $(x,y)\mapsto S_x(y)$;
it does not depend on $S$ and is again $H$-invariant. The relation
\eqref{eqhom} shows that the law $\nu_t$ of $X_t$
satisfies $\nu_t=\nu_s\ast\nu_{t-s}$. If the convolution product is
commutative on the set of $H$-invariant measures,
then $(G,H)$ is said to be a Gelfand pair.

An example of space $M$ is the hyperbolic space viewed as a subspace of
the Minkowski space, namely
\[
\hyperb^d= \bigl\{x=(x_0,x_1,
\ldots,x_d)\in\reel^{1,d};|x|=1,x_0>0 \bigr\},\qquad
|x|^2=x_0^2-\sum
_{i=1}^dx_i^2,
\]
with $o=(1,0,\ldots,0)$. It can be viewed as $G/H$, where
$G=\SO^+(1,d)$ is the restricted Lorentz group of linear
transformations of $\reel^{1,d}$ which preserve the pseudo-norm, the
time direction and the space orientation, and
$H=\{h\in G;h(o)=o\}\sim \SO(d)$. Then it is known that $(G,H)$ is a
Gelfand pair.

For our result, we need the functions
\[
\chi_M(x)=\chi_G(S_x),\qquad \am(x)= |
\Ad_{S_x} |
\]
which do not depend on $S$ because $\chi_G=1$ on $H$ and the inner
product of $\lieal$ is $\Ad(H)$-invariant.

%th5 #&#
\begin{theorem}\label{ginv}
On $M=G/H$ for $H$ compact, let $X_t$ be a $G$-invariant Markov process
without Brownian part, with $X_0=o$, the L\'{e}vy
measure $\mu_X$ of which satisfies the non degeneracy assumption
\eqref{nondlie} where $\exp^{-1}=\exp_o^{-1}$ denotes
the inverse Riemannian exponential function based at $o$, and the
additional condition \eqref{symlie} if $\alpha=1$. If
$\alpha<1$ suppose moreover that $X_t$ is a pure jump process. Then
the law of $X_t$, $t>0$, is absolutely continuous
with respect to the $G$-invariant measure $\haar^{\;M}$. Let $\ell\ge0$. If
%
%e48 #&#
\begin{equation}\label{intvcchi}
\int_{V^c}\chi_M(x)
\am(x)^j\mu_X(\mathrm{d}x)<\infty
\end{equation}
for a relatively compact neighbourhood $V$ of $o$ and for $j\le\ell$,
then the density is in $C_{b,\leftarrow}^\ell$
(see Section \ref{homog} for the definition). In particular, if
%
%e49 #&#
\begin{equation}\label{intvca}
\int_{V^c}\am(x)^j
\mu_X(\mathrm{d}x)<\infty
\end{equation}
for any $j$, the density is in $C_{b,\leftarrow}^\infty$. If $(G,H)$
is a Gelfand pair, conditions \eqref{intvcchi} or
\eqref{intvca} are not needed.
\end{theorem}

\begin{pf}
We write $X$ as $X_t=\pi(\Xi_t')$ for a left L\'{e}vy process $\Xi'$
with L\'{e}vy measure $\mu_{\Xi'}$ given by \eqref{muxia};
denote by $\gener\,'$ its infinitesimal generator. The question is to
know whether the non degeneracy condition
\eqref{nondlie} for $X$ can be translated into the similar condition
for $\Xi'$ (and similarly for \eqref{symlie} when
$\alpha=1$). Recall that $\lieal$ is written as the orthogonal sum of
$\lieh$ and $\liep$, and we can choose the
section $S$ such that $S_x\in\exp\liep$ for $x$ in a neighbourhood
of $o$; then $S$ is uniquely determined and smooth
on a maybe smaller neighbourhood; the measure $\mu_X$ on a
neighbourhood of $o$ is therefore transported to a measure
$\mu'$ on a neighbourhood of 0 in $\liep$, and $\mu'$ satisfies
Assumption~ \ref{aslevy}. If $S_x=\exp\lambda$, then
\[
\exp^{-1}\bigl(hS_xh^{-1}\bigr)=
\Ad_h\lambda+\mathrm{O}\bigl(|\lambda|^2\bigr).
\]
On the other hand, we have from \eqref{muxia} that
\begin{eqnarray*}
I(\rho)&=&\int_{\{|\exp^{-1}\xi|\le\rho\}} \bigl\langle\exp^{-1}\xi, u
\bigr\rangle^2\mu_{\Xi'}(\mathrm{d}\xi)
\\
&=&\int\int_{\{|\exp^{-1}(hS_xh^{-1})|\le\rho\}} \bigl\langle\exp^{-1}
\bigl(hS_xh^{-1}\bigr),u \bigr\rangle^2
\haar^{\;H}(\mathrm{d}h)\mu_X(\mathrm{d}x),
\end{eqnarray*}
so obtaining the lower and upper bounds \eqref{nondlie} for $I(\rho)$
is equivalent to estimating
\[
I'(\rho)=\int\int_{H\times\{|\lambda|\le\rho\}} \langle
\Ad_h\lambda,u \rangle^2\haar^{\;H}(\mathrm{d}h)
\mu'(\mathrm{d}\lambda).
\]
The upper bound follows easily since $\mu'$ satisfies Assumption \ref
{aslevy}, and for the lower bound, we can restrict
the domain of integration $H$ to a small neighbourhood of the unity on
which $|\Ad_h\lambda-\lambda|\le c|\lambda|$ for
$c$ arbitrarily small. We deduce that \eqref{nondlie} holds true for
the process $\Xi'$ but the lower bound is only for
$u$ in $\liep$. In order to obtain all of $\lieal$, we add extra
independent noise in $\Xi'$. Let $\Xi''$ be a left
L\'{e}vy process on $H$ with generator $\gener\,''$ satisfying the
conditions of Theorem \ref{llevy}; it is therefore
associated to a measure $\mu''$ on a neighbourhood of 0 in $\lieh$,
satisfying \eqref{nondlie}; we have
$\gener\,''(F\circ\pi)=0$ for any smooth $F$ on $M$. Let $\Xi$ be the
process with generator $\gener=\gener\,'+\gener\,''$.
Then $\gener(F\circ\pi)=\gener\,'(F\circ\pi)$, so $X_t=\pi(\Xi_t')$ can also be written as $X_t=\pi(\Xi_t)$, and $\Xi_t$
is a left L\'{e}vy process associated to the measure $\mu'+\mu''$
which now satisfies \eqref{nondlie} for any $u$ in
$\lieal$. Condition \eqref{symlie} is similarly extended to $\Xi$
when $\alpha=1$.

Thus, we deduce from Theorem \ref{llevy} that $\Xi_t$ has a density
$p_\leftarrow(t,\cdot)$ with respect to
$\haar_\leftarrow^{\;G}$. We have from \eqref{plr} that the density of
$X_t$ with respect to $\haar^{\;M}$ is given by
\[
p(t,x)=\int_Hp_\leftarrow(t,S_xh)
\haar^{\;H}(\mathrm{d}h).
\]
Condition \eqref{intvcchi} for $X$ implies \eqref{chigad} for $\Xi'$
because $\chi_G(h)=|\Ad_h|=1$ for $h$ in $H$, so
that $\chi_G(hS_xh^{-1})=\chi_M(x)$ and $|\Ad_{(hS_xh^{-1})}|=\am
(x)$, and \eqref{chigad} always holds for $\Xi''$ for
the same reason. Thus, Theorem \ref{llevy} can also be applied to $\Xi
$ for the smoothness of the law. In the case
$\ell=0$, the continuity of $p$ follows from the continuity of
$p_\leftarrow$ and the fact that we can choose a $S$
which is smooth in a neighbourhood of $x$. In the case $\ell=1$, we have
\[
p\bigl(t,\pi(g)\bigr)=\int_Hp_\leftarrow(t,S_{\pi(g)}h)
\haar^{\;H}(\mathrm{d}h)=\int_Hp_\leftarrow(t,gh)
\haar^{\;H}(\mathrm{d}h)
\]
from the left invariance of $\haar^{\;H}$, and we deduce from the
definition \eqref{dlfxu} that
\[
D_\leftarrow p(t,x)u=\int_HD_\leftarrow
p_\leftarrow(t,S_xh)\Ad_h^{-1}u
\haar^{\;H}(\mathrm{d}h),
\]
so
\[
\bigl|D_\leftarrow p(t,x)\bigr |\le\int_H\bigl |D_\leftarrow
p_\leftarrow(t,S_xh) \bigr|\haar^{\;H}(\mathrm{d}h).
\]
The study of higher order derivatives is similar.

If $(G,H)$ is a Gelfand pair, we write a decomposition $\mu_X=\mu_X^\flat+\mu_X^\sharp$ where $\mu_X^\flat$ is the
restriction to a $H$-invariant relatively compact neighbourhood of $o$;
this corresponds to a decomposition
$\gener=\gener^{\,\flat}+\gener^{\,\sharp}$, and $X$ can be viewed as the
process with generator $\gener^{\,\flat}$ interlaced with
big jumps described by $\gener^{\,\sharp}$. Conditionally on the times of
the $N_t$ big jumps before $t$, the law of $X_t$
is therefore the convolution of $(2N_t+1)$ $H$-invariant laws. From the
commutativity of the convolution, all the big
jumps can be put together, and we can write $X_t$ in law for $t$ fixed
as $S_{X_t^\sharp}(X_t^\flat)$. The law of
$X_t^\flat$ is in $C_{b,\leftarrow}^\infty$, and this smoothness is
preserved under the action of $S_{X_t^\sharp}$.
\end{pf}

If $(G,H)$ is a Gelfand pair, the technique of previous proof for
putting together big jumps can be extended to other
cases. For instance, we can obtain the smoothness of the law if the L\'
{e}vy measure is the sum of two measures, and only
one of them satisfies the assumptions. This implies that the upper
bound in \eqref{nondlie} can be weakened, as in
\cite{liaowang07}.

We can also consider the class of processes $X_t=\pi(\Xi_t)$, for
right L\'{e}vy processes on $G$. There are Markov
processes with semigroup $P_tf(x)=\esp f(\Xi_t(x))$. Noticing that
$\Xi_t$ can also be viewed at fixed time as the
value of a left L\'{e}vy process, we again apply left differential
calculus and immediately obtain from the above proof the
following result.

%th6 #&#
\begin{theorem}\label{rlevy}
On $M=G/H$ for $H$ compact, let $X_t=\pi(\Xi_t)=\Xi_t(o)$ for a
right L\'{e}vy process $\Xi_t$ on $G$. We suppose that the
left L\'{e}vy process having the same generator at $e$ as $\Xi$
satisfies the assumptions of Theorem \ref{llevy} for some
$\ell$. Then $X_t$ has a $C_{b,\leftarrow}^\ell$ density with
respect to the $G$-invariant measure $\haar^{\;M}$.
\end{theorem}

%s4.4 #&#
\subsection{Case 2: Semi-direct product}

We consider as in Theorem \ref{rlevy} processes $X_t=\pi(\Xi_t)=\Xi_t(o)$ where $\Xi_t$ is a right L\'{e}vy process on $G$,
but do not assume that $H$ is compact. Instead, we suppose that $G$ is
a semi-direct product as described in
\eqref{gsmh}, and we apply the right differential calculus; recall
that the right Haar measure $\haar_\rightarrow^{\;M}$ of
$M$ is relatively invariant under the action of $G$ with multiplier
$\chi$ given by \eqref{multipm}.

A typical example is when $M$ is the additive group $\reel^d$ and $G$
is the affine group. Then $H=GL(d)$.

%th7 #&#
\begin{theorem}\label{semid}
On $M=G/H$ for $G=S(M)\rtimes H$, let $X_t=\pi(\Xi_t)=\Xi_t(o)$ for
a right L\'{e}vy process $\Xi_t$ on $G$. We suppose
that the L\'{e}vy measure $\mu_\Xi$ of $\Xi$ satisfies \eqref
{nondlie}, and the additional condition \eqref{symlie} if
$\alpha=1$, and that $\Xi$ is a pure jump process if $\alpha<1$.
Assume also that
%
%e50 #&#
\begin{equation}\label{intvcm}
\int_{V^c}\frac1{\chi(g)} \bigl|
\Ad_g^{-1} \bigr|^j\mu_\Xi(\mathrm{d}g)<\infty
\end{equation}
for $j\le\ell$, where $V$ is a relatively compact neighbourhood of
$e$, and $\chi$ is given by \eqref{multipm}. Then
$X_t$ has a $C_{b,\rightarrow}^\ell$ density with respect to the
right Haar measure $\haar_\rightarrow^{\;M}$.
\end{theorem}

\begin{pf}
Let us consider $\Xi_t$. We know that it can be viewed as the solution
of an equation driven by a $\lieal$-valued L\'{e}vy
process $\Lambda$. We denote by $\Xi^\eps$ or $\Xi(\eps)$ the same
process when jumps of $\Lambda$ greater than some
$\eps$ (for some norm) have been removed, and by $X^\eps$ or $X(\eps
)$ its projection on~$M$. We have seen in
\eqref{estdjp} that we have an estimate for the density of $\Xi^\eps
$ and its derivatives at $e$ involving
$\proba[\Xi_t^\eps\in V_G]$ for a relatively compact neighbourhood
$V_G$ of $e$, uniformly in the initial condition.
From the right invariance of the process, we also have an estimate for
the density at $g$ involving
$\proba[\Xi_t^\eps\in V_Gg]$, uniform in $g$. In particular, we can write
%
%e51 #&#
\begin{equation}\label{drjp}
\bigl|D_\rightarrow^jp_\rightarrow^{\Xi(\eps)}
\bigl(t,S_xhg_0^{-1}\bigr) \bigr| \le
C_j t^{-(d+j)/\alpha}\proba\bigl[\Xi_t^\eps\in
V_GS_xhg_0^{-1}\bigr]
\end{equation}
for any $g_0\in G$, $x\in M$, $h\in H$. We want to estimate the
integral of this quantity with respect to~$h$. Let
$V_M$ be a relatively compact neighbourhood of $o$; then $V_G'=V_G
S(V_M)^{-1}$ is a relatively compact neighbourhood
of $e$. For $x$ in $M$, let $V_M^x=x^{-1}\cdot V_M\cdot x$. We have
$V_G\subset V_G'S_xS_yS_x^{-1}$ for $y\in V_M^x$,
so
\[
\proba\bigl[\Xi_t^\eps\in V_GS_xhg_0^{-1}
\bigr]\le\proba\bigl[\Xi_t^\eps\in V_G'S_xS_yhg_0^{-1}
\bigr]
\]
for any $y\in V_M^x$ and $h\in H$. By taking the mean value on $V_M^x$,
we obtain
\[
\proba\bigl[\Xi_t^\eps\in V_GS_xhg_0^{-1}
\bigr]\le\frac1{\haar_\rightarrow^{\;M}\bigl(V_M^x
\bigr)}\int_{V_M^x}\proba\bigl[\Xi_t^\eps
\in V_G'S_xS_yhg_0^{-1}
\bigr]\haar_\rightarrow^{\;M}(\mathrm{d}y),
\]
and $\haar_\rightarrow^{\;M}(V_M^x)=\haar_\rightarrow^{\;M}(V_M)/\chi_M(x)$. Thus,
%
%e52 #&#
\begin{eqnarray}\label{inthp}
&&\int_H\proba\bigl[\Xi_t^\eps
\in V_GS_xhg_0^{-1}\bigr]
\haar_\rightarrow^{\;H}(\mathrm{d}h)
\nonumber
\\
&&\quad\le C\chi_M(x)\int\int_{V_M^x\times H}\proba\bigl[
\Xi_t^\eps\in V_G'S_xS_yhg_0^{-1}
\bigr]\haar_\rightarrow^{\;M}(\mathrm{d}y)\haar_\rightarrow^{\;H}(\mathrm{d}h)
\nonumber
\\[-8pt]
\\[-8pt]
&&\quad\le C\chi_M(x)\int_G\proba\bigl[
\Xi_t^\eps\in V_G'S_xgg_0^{-1}
\bigr]\haar_\rightarrow^{\;G}(\mathrm{d}g)
\nonumber
\\
&&\quad=C\chi_M(x) \esp\bigl[\haar_\rightarrow^{\;G}
\bigl(S_x^{-1}\bigl(V_G'
\bigr)^{-1}\Xi_t^\eps g_0 \bigr)
\bigr] =C\haar_\rightarrow^{\;G}\bigl(\bigl(V_G'
\bigr)^{-1}\bigr),\nonumber
\end{eqnarray}
where we have used \eqref{fgrg} in the second inequality and $\chi_G(S_x^{-1})=1/\chi_G(S_x)=1/\chi_M(x)$ in the last
equality. The density of $\Xi_t^\eps g_0$ at $g$ is $p_\rightarrow^{\Xi(\eps)}(t,gg_0^{-1})$, and similarly for its
right invariant derivatives. The law of the variable $\pi(\Xi_t^\eps
g_0)$ is the law of $X_t^\eps$ with initial
condition $x_0=\pi(g_0)$, and \eqref{phpr} becomes
\[
p_\rightarrow^{X(\eps)}(t,x_0,x)=\int_Hp_\rightarrow^{\Xi(\eps
)}
\bigl(t,S_xhg_0^{-1}\bigr)
\haar_\rightarrow(\mathrm{d}h).
\]
By differentiating, we have
\[
\bigl|D_\rightarrow^jp_\rightarrow^{X(\eps)}(t,x_0,x)
\bigr| \le\int_H \bigl|D_\rightarrow^jp_\rightarrow^{\Xi(\eps
)}
\bigl(t,S_xhg_0^{-1}\bigr) \bigr|
\haar_\rightarrow^{\;H}(\mathrm{d}h) \le C_j t^{-(d+j)/\alpha}
\]
from \eqref{drjp} and \eqref{inthp}. Thus, the smoothness of the law
of $X^\eps$ is proved. We now have to take into
account big jumps with the technique of Lemma \ref{bjxy}, by
considering $X^\eps$ interlaced with the big jumps, and
letting $(\tau_K,\tau_{K+1})$ be the longest interval without small
jumps. Previous argument shows that, conditionally
on the times of big jumps, the variable $X_{\tau_{K+1}-}$ has a
density $p_\star$ which is in $C_{b,\rightarrow}^\ell$,
and its derivatives are of order $(\tau_{K+1}-\tau_K)^{-(d+j)/\alpha
}$. The variable $X_t$ is then obtained from the
action of $\Upsilon=\Xi_t\Xi_{\tau_{K+1}-}^{-1}$, so its density is
\[
p_\rightarrow^X(t,x_0,x)=\esp\bigl[p_\star
\bigl(\Upsilon^{-1}(x)\bigr)/\chi (\Upsilon) \bigr].
\]
The variable $1/\chi(\Upsilon)$ is conditionally integrable (given
the times of big jumps) if \eqref{intvcm} holds for
$j=0$, so the theorem can be proved for $\ell=0$ by the technique of
Lemma \ref{bjxy}; the case $\ell\ge1$ is similar
by applying \eqref{drfg} for the derivatives of $p_\star(\Upsilon^{-1}(x))$.
\end{pf}

%s5 #&#
\section{Examples}\label{exam}

We here give some examples, and also some counterexamples where the
``big jumps'' condition (Assumption \ref{asjump})
does not hold, and the smoothness of the density fails.

%s5.1 #&#
\subsection{Isotropic jumps}\label{isotrop}

We have assumed that $X$ is solution of an equation driven by some L\'
{e}vy process $\Lambda$, but usually, jumps are often
described by the L\'{e}vy kernel $\mu_x$, the image of $\mu$ by
$\lambda\mapsto a(x,\lambda)$. It is not easy to know when
some Markov process with some L\'{e}vy kernel can be represented as the
solution of an equation of our type. We have
already seen in Section \ref{jcoeff} how it is possible to deal with a
finite part of the L\'{e}vy kernel. We now give the
example of \cite{applebaumestr00} where this is globally possible.
Let $M$ be a complete Riemannian manifold, and
suppose that $\mu_x$ is the image by $(r,u)\mapsto\exp_x(ru)$ of
$\mu_R\otimes\nu_x$, where $\mu_R$ is a measure on
$(0,\infty)$ (radial part), and $\nu_x$ is the uniform probability
measure on the unit sphere of $T_xM$ (angular part);
this means that we choose a direction uniformly in the unit sphere,
then go along a geodesic in that direction, at a
distance chosen according to $\mu_R$. Such a $\mu_x$ is singular if
$\mu_R$ is singular.

In order to construct an equation for this process, as explained in
\cite{applebaumestr00}, we lift it to the bundle
$O(M)$ of orthonormal frames, as this is classically done in the
Eells--Elworthy--Malliavin construction of the Brownian
motion. Points of this bundle can be written as $\xi=(x,g)$ for $x\in
M$ and $g\dvtx \reel^d\to T_xM$ is an orthogonal
linear map; we put $\pi(\xi)=x$. Then, for $\lambda\in\reel^d$, we
can define $a(\xi,\lambda)$ for $\xi=(x,g)$ by
\[
\pi\bigl(a(\xi,\lambda)\bigr)=\exp_x(g\lambda),
\]
and the frame at $\pi(a(\xi,\lambda))$ is deduced from the frame $g$
at $x$ by parallel translation along the geodesic
$ (\exp_x(g\lambda t);0\le t\le1 )$. Let $\Xi$ be the solution of
the equation on $O(M)$ with this
coefficient $a$ and with $b=0$, driven by a symmetric L\'{e}vy process
$\Lambda$ with L\'{e}vy measure $\mu=\mu_R\otimes\nu$,
where $\nu$ is the uniform measure on the unit sphere of $\reel^d$.
Then $X=\pi(\Xi)$ is the process that we are
looking for. The process $\Xi$ can be viewed as the horizontal process
above $X$. Notice that if $(e_i)$ is the
canonical basis of $\reel^d$, the vector fields $\abar(\cdot)e_i$ are the
canonical horizontal vector fields on $O(M)$, and
the equation of $\Xi$ is a canonical equation, since $a$ is obtained
from $\abar$ by means of~\eqref{canonic}.

However, the surjectivity of $\abar(\xi)$ cannot be satisfied, since
the dimension of $O(M)$ which is $d(d+1)/2$ is
greater than the dimension $d$ of $\reel^d$. Nevertheless, we can add
extra noise in $\Xi$ without modifying the law of
$X$ and get this non degeneracy condition; more precisely, the extra
noise acts vertically on the process $\Xi$. We
enlarge the space $\reel^d$ of the L\'{e}vy process as $\reel^d\times
O(d)$, put $a(x,\lambda,e)=a(x,\lambda)$, and let
$a(x,0,g)$ be the vertical transformation which modifies the frame by
composing it with $g$. We let $\Lambda'$ be an
independent L\'{e}vy process on $O(d)$, and consider the equation
driven by $(\Lambda,\Lambda')$.

We choose the reference measure on $O(M)$ which, when projected on $M$
is the Riemannian measure, and which is on each
fibre the uniform measure (normalised measure invariant under the
action of $O(d)$). We consider the process $\Xi$ on
$O(M)$ the initial condition $\Xi_0$ of which has uniform law on the
fibre above some $x_0$. We deduce from Theorem
\ref{mainth} the smoothness of the law of $\Xi_t$ on $O(M)$ if $\mu_R$ has bounded support, and
$\int_{\{|\lambda|\le\rho\}}|\lambda|^2\mu_R(\mathrm{d}\lambda)$ is
bounded below and above by constants times $\rho^{2-\alpha}$
as $\rho\downarrow0$. Conditionally on $X_t=\pi(\Xi_t)$, the
variable $\Xi_t$ has uniform law on the fibre above $X_t$,
so the density of $\Xi_t$ is a function of $x$ only; this is also the
density of $X_t$, so $X_t$ has a smooth law too.

When $M$ is the sphere or the hyperbolic space, then $M$ is a
Riemannian symmetric space, and we are in the framework
of a $G$-invariant process on $M=G/H$, where $(G,H)$ is a Gelfand pair.
Thus, the smoothness of the law holds on
$\hyperb^d$ also in the case of unbounded jumps (Theorem~\ref{ginv}).
The hyperbolic space is diffeomorphic to
$\reel^d$; however, if we use normal coordinates, we cannot apply the
theorem of \cite{picard96} for $\reel^d$ because
the ellipticity condition is not satisfied. The matrix $(\abar\abar^\star)^{-1}$, which is the hyperbolic Riemannian
metric, explodes indeed exponentially fast at infinity, whereas at most
polynomial growth was assumed in
\cite{picard96}.

%s5.2 #&#
\subsection{Affine transformations}

Let $G$ be the group of affine transformations of $\reel^d$; it enters
the framework of \eqref{gsmh} as
$G=\reel^d\rtimes GL(d)$, where $g=(g_2,g_1)$ is the map
$g(x)=g_1x+g_2$. The Lie algebra is the set
$\reel^d\times{\mathfrak gl}(d)$, and if $u=(u_2,u_1)$, we have
\[
g\exp(\eps u)g^{-1}(x)= x+\eps\bigl(g_1u_1g_1^{-1}x+g_1
\bigl(u_2-g_1u_1 g_1^{-1}g_2
\bigr) \bigr)+\mathrm{O}\bigl(\eps^2\bigr),
\]
so
\[
\chi_G(g)= |\det g_1 |, \qquad \Ad_gu=
\bigl(g_1u_1g_1^{-1},g_1
\bigl(u_2-g_1u_1 g_1^{-1}g_2
\bigr) \bigr).
\]
We deduce from Theorem \ref{llevy} the uniform smoothness of the laws
of a left L\'{e}vy process $\Xi_t$ and its projection
$X_t=\pi(\Xi_t)$ on $\reel^d$, under the nondegeneracy condition on
the small jumps (with additional conditions if
$\alpha\le1$), and if the moments of $|\Ad_g|$ for the big jumps
part are finite. Under the similar condition on
$|\Ad_g^{-1}|$, we can also consider the right L\'{e}vy process, and
deduce the smoothness of its $\reel^d$ component from
Theorem \ref{semid}; in this case $H=GL(d)$ is unimodular, so $\chi
=\chi_G$.

Notice that $G$ is not unimodular, so without the assumption on big
jumps, we are not even sure of the local
boundedness of the density of $X_t$, except when we can apply Theorem
\ref{smoothj} (smooth L\'{e}vy measure). In order to
find a counterexample, we are going, for simplicity, to consider a
subgroup of $G$.

So let now $G$ be the group of transformations of $\reel^d$ generated
by translations and by dilations of rate $\mathrm{e}^n$,
$n\in\entier$. Thus $G=\reel^d\rtimes\entier$, where $(y,n)$
corresponds to the transformation $z\mapsto \mathrm{e}^nz+y$. The
composition law of this group is
\[
(y_2, n_2).(y_1, n_1)=
\bigl(\mathrm{e}^{n_2}y_1+y_2, n_1+n_2
\bigr).
\]
Its Lie algebra is Abelian and can be identified to $\reel^d$.
Denoting by $\mathrm{d}n$ and $\mathrm{d}y$ the counting measure on
$\entier$ and the Lebesgue measure on $\reel^d$, the Haar measures on
$G$, the modulus and the adjoint representation
are
\[
\haar_\leftarrow^{\;G}(\mathrm{d}y\, \mathrm{d}n)=\mathrm{e}^{-nd}\,\mathrm{d}y\, \mathrm{d}n,\qquad\!
\haar_\rightarrow^{\;G}(\mathrm{d}y\, \mathrm{d}n)= \mathrm{d}y\, \mathrm{d}n,\qquad\!
\chi_G(y,n)=\mathrm{e}^{nd},\qquad\!
\Ad_{(y,n)}u=\mathrm{e}^nu.
\]

We consider on $G$ left L\'{e}vy processes $X_t=(Y_t, N_t)$ such that
$N_t$ is a random walk on $\entier$, and $Y_t$ is
deduced from a L\'{e}vy process $\Lambda_t$ on $\reel^d$, independent
of $N$, by means of
\[
Y_t=\int_0^t \mathrm{e}^{N_s}\, \mathrm{d}
\Lambda_s.
\]
We want to study the density near $(0,0)$ of $X_1$ with respect to
$\haar_\leftarrow^{\;G}$, or equivalently the density
near 0 of $Y_1$ restricted to $\{N_1=0\}$. Suppose that $\Lambda$ is
an isotropic stable process. Then
\[
Y_1\sim\Lambda\biggl(\int_0^1 \mathrm{e}^{\alpha N_s}\,\mathrm{d}s
\biggr)\sim\biggl(\int_0^1\mathrm{e}^{\alpha
N_s}\,\mathrm{d}s
\biggr)^{1/\alpha}\Lambda_1,
\]
and its conditional density at 0 given $N$ is
%
%e53 #&#
\begin{equation}\label{pyon}
p^Y(0|N)=c_\Lambda\biggl(\int
_0^1 \mathrm{e}^{\alpha N_s}\,\mathrm{d}s \biggr)^{-d/\alpha},
\end{equation}
where $c_\Lambda$ is the density of $\Lambda_1$ at 0. Let $p_n$ be
the mass of the L\'{e}vy measure of $N$ at $n$, with
$\sum p_n<\infty$. Let $A_n$, $n\ge1$, be the following event: the
process $N$ has exactly two jumps before time 1, one
jump of size $-n$ followed by a jump of size $+n$. Then the probability
of $A_n$ is $c_N p_{-n}p_n$ with
$c_N=\frac12\exp-\sum p_j$; conditionally on $A_n$, the times of the
two jumps are obtained from two independent
variables uniformly distributed on $[0,1]$, and by letting the negative
jump be the smallest one, and the positive jump
the largest one. The law of $X_1$ restricted to $A_n$ is smooth; by
applying \eqref{pyon} and by denoting $c=c_\Lambda
c_N$, its density at $(0,0)$ is
\begin{eqnarray*}
p_{A_n}(0,0)&=&2c p_{-n}p_n\int
_0^1\int_0^t
\bigl(s+(t-s)\mathrm{e}^{-n\alpha}+1-t \bigr)^{-d/\alpha}\,\mathrm{d}s\, \mathrm{d}t
\\
&=&2c p_{-n}p_n\int_0^1(1-s)
\bigl(1+s\bigl(\mathrm{e}^{-n\alpha}-1\bigr) \bigr)^{-d/\alpha}\,\mathrm{d}s
\\
&\ge& c p_{-n}p_n\mathrm{e}^{-2n\alpha} \bigl(1-
\bigl(1-\mathrm{e}^{-n\alpha}\bigr)^2 \bigr)^{-d/\alpha}
\\
&\ge& c p_{-n}p_n\mathrm{e}^{n(d-2\alpha)}.
\end{eqnarray*}
The law of $X_1$ restricted to $\bigcup_{1\le k\le n}A_k$ is smooth
with density $\sum_{k=1}^np_{A_k}$, so for any
neighbourhood $U$ of $(0,0)$, the density $p$ of $X_1$ satisfies
\[
\operatorname{\ess\sup}\limits_Up\ge c\sum_{k=1}^np_{-k}p_k\mathrm{e}^{k(d-2\alpha)}.
\]
If the series diverges, then $X_1$ does not have a locally bounded
density, and examples can be constructed as soon as
$d>2\alpha$.\vadjust{\goodbreak}

Choose for instance
%
%e54 #&#
\begin{equation}\label{pnn}
p_n=\mathrm{e}^{-\beta n},\qquad p_{-n}=\mathrm{e}^{-\sigma n}
\end{equation}
for $n\ge1$, $\beta>0$, $\sigma>0$. In this case the density is not
locally bounded if $d\ge2\alpha+\beta+\sigma$. On
the other hand, Theorem \ref{llevy} implies that if $d<\beta$, the
density is bounded uniformly with respect to the
initial condition.

Consider now the law of $Y_1$. We can consider the right L\'{e}vy
process $X'$ with the same generator at $e$ as $X$; if
$N$ and $\Lambda$ are the same independent random walk and stable
isotropic processes, $X'$ can be obtained as
$X'=(Y',N)$, where $Y'$ the solution of $Y_{t+\mathrm{d}t}'=a(Y_t',\mathrm{d}N_t,\mathrm{d}\Lambda_t)$, and $a(y,n,\lambda)=\mathrm{e}^ny+\lambda$. We have
the equalities in law $X_1'\sim X_1$ and $Y_1'\sim Y_1$.

The study of the law of $Y'$ has been the subject of Theorem \ref
{semid}. For the example \eqref{pnn}, this theorem
implies that $Y_1'\sim Y_1$ has a bounded density if $d<\sigma$. On
the other hand, we have just verified that this
density is not locally bounded if $d\ge2\alpha+\beta+\sigma$;
actually, this condition can be improved because we now
study $Y_1$ instead of $X_1$. The event $A_n$ can be replaced by
$A_n'$: the process $N$ has exactly one jump of size
$-n$ before time 1, and we deduce that the density is not locally
bounded if $d\ge\alpha+\sigma$.

%s5.3 #&#
\subsection{Killed processes}\label{killed}

Consider a process which satisfies our sufficient conditions for the
smoothness of the law in some manifold, say the
affine space $\reel^d$, but which is killed at the exit from some open
subset~$M$. The existence of a locally bounded
density is preserved, but it appears that this killing can destroy the
smoothness of this density. We have seen in
\eqref{killj} that Assumption \ref{ascoeff} is preserved, but
Assumption~\ref{asjump} may fail; a problem can occur
when jumps are allowed to enter the subset (coming from the infinity of
$M$), or at least (by applying Theorem
\ref{smoothj}) if the non-smooth part of these jumps are allowed to
enter the subset. If directions of non-smooth jumps
lie in some closed cone of the vector space $\reel^d$, then smoothness
is preserved if the obstacle $M^c$ is such a
cone based at some point of $\reel^d$; in the general case however,
roughly speaking, the obstacle can produce
singularity behind it. Let us give a one-dimensional example.

Let $\Lambda_t=\Lambda_t^0-\Lambda_t^1$ where $\Lambda_t^0$ and
$\Lambda_t^1$ are real, independent, and are
respectively, a symmetric stable L\'{e}vy process with L\'{e}vy measure
$|\lambda|^{-\alpha-1}\,\mathrm{d}\lambda$ and a standard Poisson
process. Let $X_t$ be the process $\Lambda_t$ killed when it quits
$M=(-\infty,1)$. It is solution of the equation
corresponding to $b=0$ and $a(x,\lambda)=x+\lambda$ if $x+\lambda
<1$, equal to $\infty$ otherwise; the initial
condition is $x_0=0$.

Notice that the similar process without $\Lambda_t^1$ has smooth
densities from Theorem \ref{smoothj}. If $\Lambda_t^1$
is added instead of subtracted, the law should also be smooth since the
appended jumps do not enable the process to
enter $M$. We are now going to prove that in our framework, the law of
$X_1$ is not $C^1$ as soon as $\alpha\le2/3$.

Let $\tau$ be the lifetime of $X$, which is the first exit time of
$\Lambda$ from $M$. Then the density of $X_1$ is
%
%e55 #&#
\begin{equation}\label{pqux}
p(x)=q(1,x)-\esp\bigl[q(1-\tau,x-\Lambda_\tau)1_{\{\tau<1\}}
\bigr],
\end{equation}
where $q(t,\cdot)$ is the density of $\Lambda_t$. We have
\[
q(t,x)=\mathrm{e}^{-t}\sum_k
\frac{t^k}{k!}q_0(t,x+k),
\]
where the density $q_0$ of $\Lambda^0$ is $C^\infty$ on $(\reel_+\times\reel)\setminus\{(0,0)\}$, with bounded
derivatives out of any neighbourhood of $(0,0)$. Consequently, $q$ is
smooth on
$(\reel_+\times\reel)\setminus(\{0\}\times\entier_-)$, and
%
%e56 #&#
\begin{equation}\label{qktx}
 q_{(k)}(t,x)=q(t,x)-\mathrm{e}^{-t}\frac{t^k}{k!}q_0(t,x+k)
\end{equation}
is smooth on $\reel_+\times(-k-1,-k+1)$. We are going to study the
derivative of $p(x)$ as $x\uparrow0$.

We need some information on the joint law of the exit time $\tau$ and
the overshoot $\Lambda_\tau-1$. To this end, we
first check that $\tau$ is almost surely a time of jump of $\Lambda$;
this means that the process cannot creep upward,
see \cite{doney07}; a simple way of verifying this fact is to notice
that $\Lambda^0$ cannot creep both upward and
downward because it has no Brownian part; since it is moreover
symmetric, it creeps neither upward, neither downward,
and $\Lambda$ which is obtained by adding a process with finitely many
jumps satisfies the same property. Consider the
joint law of $(\tau,\Lambda_{\tau-},\Lambda_\tau)$, and denote by
$\saut$ the set of jumps of $\Lambda$; notice that
$\tau$ is a jump of $\Lambda^0$; then
\begin{eqnarray*}
\esp\bigl[f(\tau,\Lambda_{\tau-},\Lambda_\tau)1_{\{\tau<\infty\}}
\bigr] &=&\esp\sum_{t\in\saut}f(t,X_{t-},X_{t-}+
\Delta\Lambda_t)1_{\{
X_{t-}\ne\infty\}}1_{\{\Delta\Lambda_t\ge1-X_{t-}\}}
\\
&=&\esp\int_0^\infty\int_0^\infty
f(t,X_t,X_t+\lambda)1_{\{X_t\ne
\infty\}}1_{\{\lambda\ge1-X_t\}}
\lambda^{-1-\alpha}\,\mathrm{d}\lambda \,\mathrm{d}t
\nonumber
\end{eqnarray*}
from a key formula of stochastic calculus on Poisson spaces. Thus,
\[
\proba[\tau\in \mathrm{d}t,\Lambda_{\tau-}\in \mathrm{d}x,\Lambda_\tau\in \mathrm{d}y ]
=\proba[X_t\in dx ](y-x)^{-1-\alpha}\,\mathrm{d}t \,\mathrm{d}y
\]
on $(0,\infty)\times(-\infty,1)\times(1,\infty)$, and the density
of $(\tau,\Lambda_\tau)$ is
%
%e57 #&#
\begin{equation}\label{pity}
\zeta(t,y)=\esp\bigl[(y-X_t)^{-1-\alpha} \bigr]
\end{equation}
on $(0,\infty)\times(1,\infty)$; in particular $y\mapsto\zeta
(t,y)$ is $C_b^\infty$ on $[1+\eps,\infty)$ for any
$\eps>0$, uniformly in~$t$. Let $h\dvtx\reel\to[0,1]$ be a smooth
function which is 1 on $(-\infty,4/3]$ and 0 on
$[5/3,+\infty)$. From~\eqref{pqux}, we have
\begin{eqnarray*}
p(x)&=&q(1,x)-\int_1^\infty\int
_0^1 q(1-t,x-y)\zeta(t,y)\,\mathrm{d}t \,\mathrm{d}y
\\
&=&q(1,x)-\int_1^\infty\int_0^1
q(1-t,x-y) \bigl(1-h(y)\bigr)\zeta(t,y)\,\mathrm{d}t \,\mathrm{d}y-\pbar(x)
\end{eqnarray*}
with
\[
\pbar(x)=\int_1^\infty\int_0^1
q(1-t,x-y)h(y)\zeta(t,y)\,\mathrm{d}t \,\mathrm{d}y.
\]
The measure $(1-h(y))\zeta(t,y)\,\mathrm{d}y$ is $C_b^\infty$, uniformly in $t$,
so its convolution with the law of
$\Lambda_{1-t}$ is also smooth, and we obtain that $p+\pbar$ is
smooth. It is therefore sufficient to study the
regularity of $\pbar$.

The function $q(1-t,\cdot)$ is smooth out of $\entier_-$, so we can
differentiate $\pbar$ on $(-1/3,0)$ and obtain
\begin{eqnarray*}
D\pbar(x)&=&\int_1^{5/3}\int_0^1Dq(1-t,x-y)h(y)
\zeta(t,y)\,\mathrm{d}t \,\mathrm{d}y
\\
&=&\int_1^{5/3}\int_0^1(1-t)\mathrm{e}^{t-1}Dq_0(1-t,x+1-y)h(y)
\zeta(t,y)\,\mathrm{d}t \,\mathrm{d}y+\mathrm{O}(1)
\end{eqnarray*}
as $x\uparrow0$, because the rest involves the function $q_{(1)}$ of
\eqref{qktx} which is smooth on $(-2,0)$. The
self-similarity of $\Lambda^0$ enables to write
\[
D\pbar(x)=\int_1^{5/3}\int_0^1(1-t)^{1-2/\alpha}\mathrm{e}^{t-1}Dq_0
\biggl(1,\frac{x+1-y}{(1-t)^{1/\alpha}} \biggr)h(y)\zeta(t,y)\,\mathrm{d}t \,\mathrm{d}y+\mathrm{O}(1).
\]
The law of the stable variable $\Lambda_1^0$ is symmetric and
unimodal, so $Dq_0$ is nonnegative on $\reel_-$, and is
bounded below by a positive constant on some $[-C_2,-C_1]\subset
(-\infty,-2^{1/\alpha}]$. Thus, the double integral can
be bounded below by considering only the part where the fraction in
$Dq_0$ is in $[-C_2,-C_1]$, where $1/2\le t\le1$,
and where $1<y\le4/3$. With the change $s=1-t$, we obtain
%
%e58 #&#
\begin{equation}\label{dpbar}
D\pbar(x)\ge c\int_1^{4/3}\int
_{J(y)}s^{1-2/\alpha}\zeta(1-s,y)\,\mathrm{d}s \,\mathrm{d}y-C
\end{equation}
with
\[
J(y)= \bigl\{0<s\le1/2;C_1\le s^{-1/\alpha}(y-1-x)\le
C_2 \bigr\}.
\]
On the other hand, from \eqref{pity},
\[
\zeta(t,y)=\esp\bigl[(y-X_t)^{-1-\alpha} \bigr]\ge
\mathrm{e}^{-t}\esp \bigl[\bigl(y-X_t^0
\bigr)^{-1-\alpha} \bigr],
\]
where $\mathrm{e}^{-t}$ is the probability for $\Lambda^1$ to be 0 up to time
$t$, and $X^0$ is the process $\Lambda^0$ killed
at the exit of $M$. We have from \cite{chenkimsong10} that the law
of $X_t^0$ with initial condition $X_0^0=0$ is
bounded below and above by some positive constants times $(1-u)^{\alpha
/2}\,\mathrm{d}u$ for $0\le u<1$ and $1/2\le t\le1$, so
%
%e59 #&#
\begin{eqnarray}\label{zetas}
\zeta(1-s,y)&\ge& c\int_0^1(y-u)^{-\alpha-1}(1-u)^{\alpha
/2}\,\mathrm{d}u \nonumber
\\[-8pt]
\\[-8pt]
&\ge& c\int_{1-3(y-1)}^{1-(y-1)}\bigl((y-1)+(1-u)
\bigr)^{-\alpha-1}(1-u)^{\alpha/2}\,\mathrm{d}u \ge c'(y-1)^{-\alpha/2}\nonumber
\end{eqnarray}
for $1<y\le4/3$ and $0\le s\le1/2$. We also have
\[
J(y)= \bigl\{s>0;C_1\le s^{-1/\alpha}(y-1-x)\le C_2
\bigr\}
\]
for $1<y\le5/3$ and $-1/3\le x<0$, because $0\le y-1-x\le1$ and
$C_2^{-\alpha}<C_1^{-\alpha}\le1/2$; thus
%
%e60 #&#
\begin{equation}\label{intjy}
\int_{J(y)}s^{1-2/\alpha}\,\mathrm{d}s=C(y-1-x)^{2\alpha-2}.
\end{equation}
It follows from \eqref{dpbar}, \eqref{zetas} and \eqref{intjy} that
\[
D\pbar(x)\ge c\int_1^{4/3}(y-1-x)^{2\alpha-2}(y-1)^{-\alpha/2}\,\mathrm{d}y-C.
\]
If $\alpha\le2/3$, we obtain that $D\pbar(x)$ tends to $+\infty$ as
$x\uparrow0$, so $Dp(x)$ tends to $-\infty$ and the
law of $X_1$ is not $C^1$.

\section*{Acknowledgement}
Jean Picard is supported in part by ANR Project ProbaGeo ANR-09-BLAN-0364-01.

%suskaldyti doi

% imsref loaded by aiste.veprauskaite, 2012-06-19 09:36:14

\printhistory

\end{document}